\documentclass{amsart}

\usepackage{amssymb}

\newtheorem{theorem}{Theorem}[section]
\newtheorem{lemma}[theorem]{Lemma}
\newtheorem{corollary}[theorem]{Corollary}
\newtheorem{proposition}[theorem]{Proposition}
\newtheorem{claim}[theorem]{Claim}

\newtheorem{conjecture}[theorem]{Conjecture}

\numberwithin{equation}{section}
\newcommand{\Ai}{\text{Ai\,}}

\newcommand{\im}{\text{Im\,}}

\begin{document}
\setcounter{page}{1}

\title[Discrete polynuclear growth and determinantal processes]
{Discrete polynuclear growth and determinantal processes}
\author[K.~Johansson]{Kurt Johansson}

\address{
Department of Mathematics,
Royal Institute of Technology,
S-100 44 Stockholm, Sweden}

\email{kurtj@math.kth.se}

\begin{abstract}
We consider a discrete polynuclear growth  (PNG) process and prove a
functional limit theorem for its convergence to the Airy process. This
generalizes previous results by Pr\"ahofer and Spohn. The result
enables us to express the $F_1$ GOE Tracy-Widom distribution in terms
of the Airy process. We also show some results, and give a conjecture,
about the transversal fluctuations in a point to line last passage
percolation problem.
\end{abstract}

\maketitle

\section{Introduction and results}
\subsection{Discrete polynuclear growth}
Recently there has been interesting developments concerning certain special 
$1+1$ dimensional local random growth models. This development has its 
starting point in the new results on the longest increasing subsequence
in a random permutation, \cite{BDJ}. We will not review all these 
developments here. In this paper we consider a certain discrete growth model
called the discrete polynuclear growth (PNG) model, \cite{KrSp}, a special
version of which is closely related to the last-passage percolation 
problem studied in \cite{Jo1}. It is a discrete version of the PNG model
studied by Pr\"ahofer and Spohn, \cite{PrSp}, which can be obtained as a 
special limiting case. In the paper we will extend the results in \cite{PrSp}
to the present model and prove a stronger convergence result. We also
obtain some preliminary results on the transversal fluctuations in the
point to line version of the last-passage percolation problem, which should
have many similarities with the corresponding problems for first-passage 
percolation and directed polmers.

The {\it discrete polynuclear growth} (PNG) model is a local random
growth model defined by
\begin{equation}\label{0.16}
h(x,t+1)=\max (h(x-1,t),h(x,t),h(x+1,t))+\omega(x,t+1),
\end{equation}
$x\in\mathbb{Z}$, $t\in\mathbb{N}$, $h(x,0)=0$, $x\in\mathbb{Z}$. 
Here $\omega(x,t)$, $(x,t)\in
\mathbb{Z}\times\mathbb{N}$, are independent random variables, see
\cite{KrSp}. Typically they could be Bernoulli random variables. 
We should think of $h(x,t)$ as the {\it height} above $x$ at time $t$,
so $x\to h(x,t)$ gives an interface developing in time.
We will
treat a special case where $\omega(x,t)=0$ if $t-x$ is even or if 
$|x|>t$, and 
\begin{equation}\label{0.17}
w(i,j)=\omega(i-j,i+j-1),
\end{equation}
$(i,j)\in\mathbb{Z}_+^2$, are independent geometric random variables
with parameter $a_ib_j$,
\begin{equation}\label{0.18}
\mathbb{P}[w(i,j)=m]=(1-a_ib_j)(a_ib_j)^m,
\end{equation}
$m\ge 0$. We will mainly consider the case when $a_i=b_i=\sqrt{q}$,
$0<q<1$, $i\ge 1$, and we we do this in the rest of this section. If
we define 
\begin{equation}\label{0.19}
G(i,j)=h(i-j,i+j-1),
\end{equation}
$(i,j)\in\mathbb{Z}_+^2$, it follows from (\ref{0.16}) that
\begin{equation}\label{0.20}
G(i,j)=\max(G(i-1,j),G(i,j-1))+w(i,j),
\end{equation}
see proposition \ref{P3.01}. This leads immediately to a different formula for
$G(M,N)$, \cite{Jo1},
\begin{equation}
G(M,N)=\max_\pi \sum_{(i,j)\in\pi} w(i,j),
\notag
\end{equation}
where the maximum is taken over all up/right paths from $(1,1)$ to
$(M,N)$. We can think of $G(M,N)$ as a {\it point to point
last-passage time}. It is also natural, from the point of view of
directed polymers for example, to consider the {\it point to line
last-passage time},
\begin{equation}\label{0.21}
G_{pl}(N)=\max_{|K|<N} G(N+K,N-K).
\end{equation}
This makes it reasonable to study the process $K\to G(N+K,N-K)$,
$-N<K<N$, which, by (\ref{0.19}), is the same as $K\to h(2K, 2N-1)$,
i.e. the height curve at even sites at time $2N-1$.

Let $F_1$ and $F_2$ denote the GOE respectively GUE Tracy-Widom
largest eigenvalue distributions, \cite{TW1}. It is known, \cite{Jo1},
that there are constants $a=2\sqrt{q}(1-\sqrt{q})^{-1}$ and 
$d$ given by (\ref{0.23}) below, such that
$\mathbb{P}[G(N,N)\le aN+dN^{1/3}\xi]\to F_2(\xi)$ as $N\to\infty$,
and, \cite{BR}, $\mathbb{P}[G_{pl}(N)\le aN+dN^{1/3}\xi]\to F_1(\xi)$
as $N\to\infty$. Also, if the maximum in (\ref{0.21}) is assumed at
some point $K_N$, which need not be unique, we expect $K_N$ to be of
order $N^{2/3}$, i.e. the {\it transversal fluctuations} are of order
$N^{2/3}$. This can be seen heuristically, \cite{KrSp}, and there are
some rigorous results for a related question, \cite{Jo5},
\cite{BDMMZ}, \cite{Wi}. 
These scales motivates the introduction of a rescaled
process $t\to H_N(t)$, $t\in\mathbb{R}$, defined by
\begin{equation}\label{0.22}
G(N+u,N-u)=\frac{2\sqrt{q}}{1-\sqrt{q}}N+
dN^{1/3}H_N(\frac{1-\sqrt{q}}{1+ \sqrt{q}}duN^{-2/3}),
\end{equation}
and linear interpolation, $|u|<N$, compare with \cite{PrSp}. This is
our rescaled discrete PNG process. The constant $d$ is given by
\begin{equation}\label{0.23}
d=\frac{(\sqrt{q})^{1/3}(1+\sqrt{q})^{1/3}}{1-q}.
\end{equation}
In the limit when $q$ is small and $N$ is large, we can obtain the
continuous PNG process studied by Pr\"ahofer and Spohn,
\cite{PrSp}. We want to extend their results to the present discrete
setting and also prove a stronger form of convergence to the limiting
process, a functional limit theorem. Before we can state the theorem
we must define the limiting process which is the Airy process
introduced by Pr\"ahofer and Spohn, \cite{PrSp}.

We will approach $H_N$ by considering it as the top curve in a multilayer
PNG process, compare \cite{Jo4} and
\cite{PrSp}. This will lead to measures of the form introduced in sect. 1.2
and we
will be able to use the formulas for the correlation functions derived
there. The same methods can also be applied to Dyson's Brownian
motion, compare with \cite{FNH},
which can be obtained from $N$ non-intersecting Brownian motions. The
appropriately rescaled limit as $N\to\infty$ of the top path in
Dyson's Brownian motion converges to the Airy process, see below.
This gives some
intuition about what it looks like. Its precise definition is more
technical.

The {\it extended Airy kernel}, \cite{FNH}, \cite{Ma}, \cite{PrSp}, is
defined by
\begin{equation}\label{0.24}
A(\tau,\xi;\tau',\xi')=
\begin{cases}
\int_0^\infty e^{-\lambda(\tau-\tau')}\Ai(\xi+\lambda)\Ai(\xi'+\lambda)
d\lambda, &\text{if $\tau\ge\tau'$}\\
\int_{-\infty}^0 e^{-\lambda(\tau-\tau')}\Ai(\xi+\lambda)\Ai(\xi'+\lambda)
d\lambda,  &\text{if $\tau<\tau'$.}
\end{cases}
\end{equation}
where $\Ai(\cdot)$ is the Airy function. 
When $\tau=\tau'$ the extended
Airy kernel reduces to the ordinary Airy kernel, \cite{TW1}.

We define the {\it Airy process} $t\to A(t)$ by giving its
finite-dimensional distribution functions. Given
$\xi_1,\dots,\xi_m\in\mathbb{R}$ and $\tau_1<\dots<\tau_m$ in
$\mathbb{R}$ we define $f$ on
$\{\tau_1,\dots,\tau_m\}\times\mathbb{R}$ by
\begin{equation}
f(\tau_j,x)=-\chi_{(\xi_j,\infty)}(x).
\notag
\end{equation}
Then,
\begin{equation}\label{0.28}
\mathbb{P}[A(\tau_1)\le\xi_1,\dots,A(\tau_m)\le\xi_m] = 
\det(I+f^{1/2}Af^{1/2})_{L^2(\{\tau_1,\dots,\tau_m\}\times\mathbb{R})},
\end{equation}
where we have counting measure on $\{\tau_1,\dots,\tau_m\}$ and
Lebesgue measure on $\mathbb{R}$. 
The Fredholm determinant can be 
defined via its Fredholm expansion, see sect. 2.1
below.
We will prove in section 2.2 that $f^{1/2}Af^{1/2}$ is a trace class
operator on $L^2(\{\tau_1,\dots,\tau_m\}\times\mathbb{R})$, so this is
also a Fredholm determinant in the sense of determinants for trace class
operators.
Note that in particular
\begin{equation}\label{0.29}
\mathbb{P}[A(\tau)\le\xi]=F_2(\xi).
\end{equation}

This defines
the Airy process. It is proved in \cite{PrSp} that it has a version
with continuous paths, which also follows from the results
below.
As mentioned above, 
another way of understanding the Airy process is as follows. Let
$\lambda(t)=(\lambda_1(t),\dots,\lambda_N(t))$ with
$\lambda_1(t)<\dots<\lambda_N(t)$, be the eigenvalues in 
Dyson's Brownian motion model, \cite{Dy}, for GUE with stationary
distribution $Z_N^{-1}\Delta_N(\lambda)^2\prod_{j=1}^N\exp(-\lambda_j^2)$.
Then,
\begin{equation}
\lim_{N\to\infty}\sqrt{2}N^{1/6}(\lambda_N(N^{-1/3}t)-\sqrt{2N})=A(t)
\notag
\end{equation}
say in the sense of convergense of finite-dimensional
distributions. This can be proved using the methods of the present paper,
and using techniques from \cite{Jox}, it is possible to get an integral
formula for the (extended) correlation kernel. The details will not 
be given here. This scaling limit has been studied before, see
\cite{FNH} and references therein.

In analogy with the results of \cite{PrSp}, we can show that
the rescaled height process $H_N$ converges in finite dimensional 
distributions to the Airy process.

\begin{theorem}\label{T4.3}
Let $H_N$ be the process defined by (\ref{0.22}). Then for any fixed
$t_1,\dots,t_m$ and $\xi_1,\dots,\xi_m$,
\begin{align}\label{4.28}
&\lim_{N\to\infty} \mathbb{P}[H_N(t_1)\le\xi_1,\dots, H_N(t_m)\le
\xi_m]\\ 
&\mathbb{P}[A(t_1)\le\xi_1+t_1^2,\dots, A(t_m)\le
\xi_m+t_m^2],
\notag
\end{align}
where $A$ is the Airy process.
\end{theorem}

This result can be sharpened to a functional limit theorem.

\begin{theorem}\label{T0.4} Let $A(t)$ be the Airy process defined by
  its finite-dimensional distributions, (\ref{0.28}). Also, let
  $H_N(t)$ be defined by (\ref{0.22}) and linear interpolation. Fix
  $T>0$ arbitrary. There is a continuous version of $A(t)$ and
\begin{equation}\label{0.30}
H_N(t)\to A(t)-t^2,
\end{equation}
as $N\to\infty$ in the $\text{weak}^\ast$-topology of probability measures on
$C(-T,T)$.
\end{theorem}
The theorem will be proved in section 5.2. As a corollary to this
theorem and the results of Baik and Rains, \cite{BR}, we obtain the
following result which expresses the GOE largest eigenvalue
distribution $F_1$ in terms of the Airy process.

\begin{corollary}\label{C0.5} For all $\xi\in\mathbb{R}$,
\begin{equation}\label{0.31}
F_1(\xi)=\mathbb{P}[\sup_t(A(t)-t^2)\le\xi].
\end{equation}
\end{corollary}

The proof of (\ref{0.31}) is very indirect. It would be interesting to
see a more straightforward approach.

As discussed above we are also interested in the transversal
fluctuations of the endpoint of a maximal path in the point to line
case. In our discrete model this is not well-defined, there could be
several maximal paths. Consider the random variable
\begin{equation}\label{0.32}
K_N=\inf\{u\,;\,\sup_{t\le u} H_N(t)=\sup_{t\in\mathbb{R}} H_N(t)\},
\end{equation}
the first point that gives the maximum. The corresponding quantity for
the limiting process $H(t)=A(t)-t^2$ is
\begin{equation}\label{0.33}
K=\inf\{u\,;\,\sup_{t\le u} H(t)=\sup_{t\in\mathbb{R}} H(t)\}.
\end{equation}
We would like to show that $K_N$ converges to $K$ so that we could
call the law of $K$ the asymptotic law of transversal
fluctuations. Unfortunately we can only prove this under a very
plausible 
assumption
on the Airy process. We can show,

\begin{proposition}\label{P0.6} The sequence of random variables
  $\{K_N\}_{N\ge 1}$ is tight, i.e. given $\epsilon>0$ there is a
  $T>0$ and an $N_0$ such that
\begin{equation}
\mathbb{P}[|K_N|>T]<\epsilon
\notag
\end{equation}
for all $N\ge N_0$.
\end{proposition}

The assumption we need to make on the Airy process can be formulated
as follows.

\begin{conjecture}\label{Con0.7}
Let $H(t)=A(t)-t^2$. Then, for each $T>0$, $H(t)$ has a unique point
of maximum in $[-T,T]$ almost surely.
\end{conjecture}
If we accept this we can prove

\begin{theorem}\label{T0.8}
Assume that conjecture \ref{Con0.7} is true. Then $K_N\to K$ in
distribution as $N\to\infty$
\end{theorem}
The law of $K$ is thus a natural candidate for the law of the
transversal fluctuations. It would be interesting to find a different,
more explicit, formula for this law. Assuming the truth of the same
conjecture it may also be possible to prove that the endpoints of
all maximal paths, or asymptotically maximal paths, converge to the
same limit $K$.

By using the limit results of \cite{Bar}, proposition \ref{P3.2}  
and theorem \ref{T3.4} we can obtain the correlation functions of the
eigenvalues of the succesive minors $H^{(k)}=(h_{ij})_{1\le i,j\le k}$,
  $1\le k\le N$, of an $N\times N$ GUE matrix $H=(h_{ij})_{1\le i,j\le
    N}$. In this way it is possible to get the Airy process as an
  appropriate limit of the succesive largest eigenvalues of
  $H^{(k)}$. More details will be given in future work.

We could also get the Airy process by looking at the largest
eigenvalues of coupled GUE-matrices, which is similar to looking at
Dyson's Brownian motion model for GUE. In \cite{PrSp} Pr\"ahofer and
Spohn raised the problem of finding differential equations for
probabilities of the form (\ref{0.28}) generalizing the Painlev\'e II
formulas for (\ref{0.29}). In \cite{AdvMo} the spectrum of coupled
random matrices is studied and it would be interesting to see if the
results of this paper shed some light on this problem.

\subsection{Measures defined by products of determinants}
Probability measures given by products of determinants has been
studied in several papers, e.g. by Eynard and Mehta, \cite{EyMe}, in
connection with eigenvalue correlations in chains of matrices, by 
Forrester, Nagao and Honner, \cite{FNH}, 
in connection with Dyson's Brownian motion
model and by Okounkov and Reshetikhin, \cite{OkRe}, when introducing
the so called Schur process. The problem is to compute the correlation
functions and to show that these are given by determinants so that we
obtain a determinantal point process, \cite{So2}. The same type of
correlation functions are also obtained by Pr\"ahofer and Spohn,
\cite{PrSp}, in a cascade of continuous polynuclear growth (PNG)
models. 
We will study a class of measures which include all the above as special
cases and show that we obtain determinantal correlation functions.
As an example of the result we will in sect. 2.3 investigate random 
walks on the discrete circle using the same strategy. This will lead
to an extended discrete sine kernel, compare with \cite{OkRe}.
We will see in sect. 3 that our main topic 
the discrete PNG problem fits nicely into
this framework. This particular application is very close to the 
Schur process in \cite{OkRe}, and their results could also have been used. 
In fact, we rederive their main formulas.

For $r\in\mathbb{Z}$ let $x^r=(x_1^r,\dots,x_n^r)\in\mathbb{R}^n$ and
$\bar{x}= (x^{-M+1},\dots,x^{M-1})$, $M\ge 1$. We think of $\bar{x}$
as a {\it point configuration} in $\{-M+1,\dots,
M-1\}\times\mathbb{R}$, and we also specify fixed {\it initial}
$x^{-M}$ and {\it final} $x^M$ positions. Let
$\phi_{r,r+1}:\mathbb{R}^2\to\mathbb{C}$, $r\in\mathbb{Z}$, be given
{\it transition weights}. The {\it weight} of the configuration
$\bar{x}$ is then
\begin{equation}\label{0.1}
w_{n,M}(\bar{x})=\prod_{r=-M}^{M-1}\det
(\phi_{r,r+1}(x_i^r,x_j^{r+1}))_{i,j=1}^n. 
\end{equation}
Let $d\mu$ be a given {\it reference measure} on $\mathbb{R}$,
typically Lebesgue measure or counting measure. 
We assume that $|\phi_{r,r+1}(x,y)|\le c_r(x)d_r(y)$, where
$c_r\in L^1(\mathbb{R},\mu)$ and $d_r\in L^\infty(\mathbb{R},\mu)$,
$-M\le r<M$. This assumption is not necessary but is convenient and 
suffices for the convergence of all the objects we will encounter.
The {\it partition
  function} is
\begin{equation}\label{0.2}
Z_{n,M}=\frac
1{(n!)^{2M-1}}\int_{(\mathbb{R}^n)^{2M-1}}w_{n,M}(\bar{x})d\mu
(\bar{x}),
\end{equation}
where $d\mu(\bar{x})=\prod_{r=-M+1}^{M-1}\prod_{j=1}^nd\mu(x_j^r)$.
We will assume that $Z_{n,M}\neq 0$ so that we can define the {\it
  normalized weight}
\begin{equation}\label{0.3}
p_{n,M}(\bar{x})=\frac 1{(n!)^{2M-1}Z_{n,M}}w_{n,M}(\bar{x}).
\end{equation}
If $w_{n,M}(\bar{x})\ge 0$, this is a probability density on
$(\mathbb{R}^n)^{2M-1}$ with respect to the reference measure
$d\mu(\bar{x})$.
The $(k_{-M+1},\dots,k_{M-1})$-{\it correlation function} can now be
defined in a standard way by 
\begin{align}\label{0.4}
&R_{k_{-M+1},\dots,k_{M-1}}(x_1^{-M+1},\dots,x_{k_{-M+1}}^{-M+1},
\dots,x_1^{M-1},\dots,x_{k_{M-1}}^{M-1})\\&=
\int_{\mathbb{R}^{n(2M-1)-k}} p_{n,M}(\bar{x})\prod_{r=-M+1}^{M-1}
\frac{n!}{(n-k_r)!} \prod_{j=k_r+1}^nd\mu(x_j^r),
\notag
\end{align}
where $k=k_{-M+1}+\dots+k_{M-1}$, $0\le k_j\le n$.

Given two transition functions we define their {\it convolution} by
\begin{equation}
\phi\ast\psi (x,y)=\int_{\mathbb{R}}\phi (x,z)\psi(z,y)d\mu(z).
\notag
\end{equation}
Set
\begin{equation}
\phi_{r,s}(x,y)=(\phi_{r,r+1}\ast\cdots\ast\phi_{s-1,s})(x,y)
\notag
\end{equation}
if $r<s$ and $\phi_{r,s}\equiv 0$ if $r\ge s$.
Let $A=(A_{ij})$ be the $n\times n$ matrix with elements
$A_{ij}=\phi_{-M,M}(x_i^{-M},x_j^M)$, $1\le i,j\le n$. By repeated use
of the Heine identity:
\begin{equation}\label{Heine}
\frac
1{n!}\int_{\mathbb{R}^n}\det(\phi_i(x_j))_{i,j=1}^n
\det(\psi_i(x_j))_{i,j=1}^n d\mu(x)=
\det\left(\int_{\mathbb{R}}\phi_i(t)\psi_j(t)d\mu(t)\right)_{i,j=1}^n,
\end{equation}
we see that
$Z_{n,M}=\det A$. Hence $\det A\neq 0$ by our assumption. Define
a kernel $K^{n,M}:(\{-M+1,\dots,M-1\})\times\mathbb{R})^2\to
\mathbb{C}$ by
\begin{equation}\label{0.5}
K^{n,M}(r,x;s,y)=\tilde{K}^{n,M}(r,x;s,y)-\phi_{r,s}(x,y),
\end{equation}
where
\begin{equation}\label{0.6}
\tilde{K}^{n,M}(r,x;s,y)=\sum_{i,j=1}^n\phi_{r,M}(x,x_i^M)(A^{-1})_{ij}
\phi_{-M,s}(x_j^{-M},y). 
\end{equation}
In the case $M=1$ the kernel $\tilde{K}$ has appeared before, see
\cite{TW2},
\cite{Bor} and also \cite{Jo3}.

\begin{theorem}\label{T0.1}
The correlation functions defined by (\ref{0.4}) are given by
\begin{align}\label{0.7}
&R_{k_{-M+1},\dots,k_{M-1}}(x_1^{-M+1},\dots,x_{k_{-M+1}}^{-M+1},
\dots,x_1^{M-1},\dots,x_{k_{M-1}}^{M-1})\\&=
\det(K^{n,M}(r,x_{i_r}^r;s,x_{j_s}^s))_
{-M<r,s<M, 0\le i_r\le k_r, 0\le j_s\le k_s}
\notag
\end{align}
\end{theorem}
The determinant in the right hand side  of (\ref{0.7}) has a block
structure with the blocks given by $r,s$ and having size $k_r\times
k_s$. The theorem will be proved in section 2.1. 

A case of particular interest is when the transition weights are given
by Fourier coefficients. We are then in a situation similar to that in
\cite{OkRe}. Let $f_r(e^{i\theta})$ be a function in $L^1(\mathbb{T})$
with Fourier coefficients $\hat{f}_r$. Assume that the transition
weights are given by
\begin{equation}\label{0.9}
\phi_{r,r+1}(x,y)=\hat{f}_r(y-x),
\end{equation}
$-M\le r<M$, $x,y\in\mathbb{Z}$ and that the initial and final
configurations are given by $x_J^{-M}=x_j^M=1-j$, $j=1,\dots,n$. If we
set 
\begin{equation}
f_{r,s}(z)=\prod_{\ell=r}^{s-1} f_\ell(z),
\notag
\end{equation}
$z=e^{i\theta}$, then
\begin{equation}\label{0.10}
\phi_{r,s}(x,y)=\hat{f}_{r,s}(y-x)
\end{equation}
for $r<s$. The matrix $A$ defined above is then a Toeplitz matrix with
symbol
\begin{equation}\label{0.11}
a(z)=f_{-M,M}(z)=\prod_{\ell=-M}^{M-1} f_\ell(z).
\end{equation}
Define
\begin{equation}\label{0.12}
\tilde{\mathcal{K}}^{n,M}_{r,s}(z,w)=\sum_{x,y\in\mathbb{Z}}
\tilde{K}^{n,M}(r,x;s,y)z^xw^{-y},
\end{equation}
where $\tilde{K}^{n,M}$ is given by (\ref{0.6}). When the transition
functions and the initial and final configurations are given in this
way we are able to give a formula for the limit of this generating
function as $n\to\infty$.

\begin{proposition}\label{P0.3} Assume that $f_r(z)$ has winding number
  zero, a Wiener-Hopf factorization $f_r(z)=f^+_r(z)f^-_r(z)$ and is
  analytic in $1-\epsilon_r<|z|<1+\epsilon_r$ for some
  $\epsilon_r>0$. Furthermore, suppose that
\begin{equation}
\sum_{n\in\mathbb{Z}} |n|^\alpha |\hat{a}_n|<\infty,
\notag
\end{equation}
for some $\alpha>0$, where $\hat{a}_n$ are the Fourier coefficients of
the symbol $a(z)$ given by (\ref{0.11}). Set $\epsilon
=\min\epsilon_r$ and 
\begin{equation}\label{0.12'}
\tilde{\mathcal{K}}^M_{r,s}(z,w)=\frac z{z-w} G(z,w),
\end{equation}
where
\begin{equation}\label{0.13}
G(z,w)=\frac{\prod_{t=r}^{M-1}f^-_t(\frac 1z)
  \prod_{t=-M}^{s-1}f^+_t(\frac 1w)}{\prod_{t=-M}^{r-1}f^+_t(\frac 1z) 
\prod_{t=s}^{M-1}f^-_t(\frac 1w)}.
\end{equation}
Then, for $1-\epsilon<|w|<1<|z|<1+\epsilon$,
\begin{equation}\label{0.14}
\left|
  \tilde{\mathcal{K}}^{n,M}_{r,s}(z,w)-\tilde{\mathcal{K}}^M_{r,s}(z,w)\right| 
\le\frac{|f_{r,M}(\frac 1z)||f_{-M,s}(\frac 1w)|}{
  (|z|-1)(1-|w|)}\left( \frac 1{n^\alpha}+|w|^{n/2}+\frac
  1{|z|^{n/2}}\right). 
\end{equation}
Furthermore,
\begin{equation}\label{0.15}
\phi_{r,s}(x,y)=\frac 1{2\pi}\int_{-\pi}^\pi e^{i(y-x)\theta}
G(e^{i\theta},e^{i\theta})d\theta,
\end{equation}
for $r<s$.
\end{proposition}
The same type of formula for the limiting kernel was obtained in
\cite{OkRe}. The formula will be proved in section 2.1. This
proposition makes it possible to compute the asymptotics of the kernel
given by (\ref{0.5}) in certain cases, since it gives an integral
formula for the $n\to\infty$ limit of $K^{n,M}$.

The outline of the paper is as follows. In sect. 2 we will give some
general results for measures of the form (\ref{0.1}) and then as an
example discuss nonintersecting random walks on the discrete
circle. The next section applies the general theory to the discrete
PNG model and gives more explicit formulas. In sext. 4 asymptotic
results for the correlation kernel appearing in the PNG model are
stated and proved, and these are then applied in sect. 5 to prove the
necessary estimates needed for the functional limit theorem and the
transversal fluctuations.

\section{Determinantal measures}
\subsection{General theory}
In this section we will prove the results of section 1.2. We will prove 
theorem \ref{T0.1} using a generalization of the method of \cite{TW2}, 
\cite{Bor} for $\beta=2$ random matrix ensembles, see also \cite{Jo3}.
It is also possible to generalize the approach of \cite{EyMe}, which is 
closer to the original Dyson approach.

Let $\Lambda_M=\{-M+1,\dots,M-1\}\times\mathbb{R}$, $\lambda$ be
the counting measure on $\{-M+1,\dots,M-1\}$ and $\nu=\lambda\otimes
\mu$. Furthermore, we let $g:\Lambda_M\to\mathbb{C}$ be a bounded function
and define
\begin{equation}
Z_{n,M}[g]=\frac 1{(n!)^{2M-1}}\int_{(\mathbb{R}^n)^{2M-1}} 
\prod_{|r|< M}\prod_{j=1}^n (1+g(r,x_j^r))w_{n,M}(\bar{x})d\mu(\bar{x}).
\notag
\end{equation}
We want to compute $Z_{n,M}[g]/Z_{n,M}[0]$. Using the Heine identity
(\ref{Heine}) repeatedly we see that
\begin{align}
Z_{n,M}[g]&=\frac 1{(n!)^{2M-1}}\int_{(\mathbb{R}^n)^{2M-1}}
\det (\phi_{-M,-M+1}(x_i^{-M},x_j^{-M+1}))_{1\le i,j\le n}\notag\\
&\times\prod_{r=-M+1}^{M-1}\det((1-g(r,x_i^r))\phi_{r,r+1}(x_i^r,x_j^{r+1}))_
{1\le i,j\le n}d\mu(\bar{x})\notag\\
&=\det\left(\int_{\mathbb{R}^{2M-1}}\phi_{-M,-M+1}(x_i^{-M},t_{-M+1})
\prod_{|r|<M}(1+g(r,t_r))\right.\notag\\ 
&\left.\times\left(\prod_{r=-M+1}^{M-2}\phi_{r,r+1}(t_r,t_{r+1})
\right) \phi_{M-1,M}(t_{M-1},x_j^M)d^{2M-1}\mu(t)\right)_{1\le i,j\le n}.
\notag
\end{align}
Now,
\begin{equation}
\prod_{|r|<M}(1+g(r,t_r))=1+\sum_{\ell=1}^{2M-1}\sum_{-M<r_1<\dots<r_\ell<M}
g(r_1,t_{r_1})\dots g(r_\ell,t_{r_\ell}),
\notag
\end{equation}
and hence
\begin{align}\label{2.01}
&Z_{n,M}[g]=\det \left(A_{ij}+\sum_{\ell=1}^{2M-1}
\sum_{-M<r_1<\dots<r_\ell<M}
\int_{\mathbb{R}^\ell}\phi_{-M,r_1}(x_i^{-M},t_1)\right.\\
&\left.\left( \prod_{s=1}^{\ell-1} g(r_s,t_s)\phi_{r_s,r_{s+1}}(t_s,t_{s+1})
\right) g(r_\ell,t_\ell)\phi_{r_\ell,M}(t_\ell,x_j^M)d^\ell\mu(t)\right)_
{1\le i,j\le n},
\notag
\end{align}
where we have used the notation of sect. 1.2. If we set $g=0$ 
we obtain $Z_{n,M}[0]=Z_{n,M}=\det A$ as before. By definition 
$\phi_{r,s}=0$ if $r\ge s$, and hence we can remove the ordering of the 
$r_i$'s in (\ref{2.01}). We find,
\begin{align}\label{2.02}
&\frac{Z_{n,M}[g]}{Z_{n,M}[0]} =\det \left(\delta_{ij}+
\sum_{k=1}^n(A^{-1})_{ik}\sum_{\ell=1}^{2M-1}
\sum_{-M<r_m<M}
\int_{\mathbb{R}^\ell}\phi_{-M,r_1}(x_k^{-M},t_1)\right.\\
&\left.\left( \prod_{s=1}^{\ell-1} g(r_s,t_s)\phi_{r_s,r_{s+1}}(t_s,t_{s+1})
\right) g(r_\ell,t_\ell)\phi_{r_\ell,M}(t_\ell,x_j^M)d^\ell\mu(t)\right)_
{1\le i,j\le n}.
\notag
\end{align}
Write $\psi(u,t;v,s)=g(u,t)\phi_{u,v}(t,s)$, and define $\psi^{\ast 0}
(u,t;v,s)=\delta_{uv}\delta(t-s)$, $\psi^{\ast 1}=\psi$ and
\begin{equation}
\psi^{\ast (r+1)}(u,t;v,s)=\int_{\Lambda_M^{r}} \psi(u,t;m_1,\xi_1)
\psi(m_1,\xi_1;m_2,\xi_2)\dots \psi(m_r,\xi_r;v,s)d^r\nu(m,\xi)
\notag
\end{equation}
for $r\ge 1$. Note that, since $\phi_{r,s}=0$ if $r\ge s$, we have
$\psi^{\ast\ell}=0$ if $\ell\ge 2M-1$. This follows immediately 
from the definition. The formula (\ref{2.02}) can now be written 
\begin{align}\label{2.03}
\frac{Z_{n,M}[g]}{Z_{n,M}[0]} &=\det \left(\delta_{ij}+
\sum_{k=1}^n(A^{-1})_{ik}\int_{\Lambda_m}d\nu(u,\xi)\int_{\Lambda_m}
d\nu(v,\eta)\right.\\
&\left.\phi_{-M,u}(x_k^{-M},\xi)\left(\sum_{\ell=1}^{2M-1}\psi^{\ast(\ell-1)}
(u,\xi;v,\eta)\right)g(v,\eta)\phi_{v,M}(\eta,x_j^M)\right)_{i,j=1,\dots,n}.
\notag
\end{align}

If $K(x,y)$ is an integral kernel on $L^2(\Omega,\mu)$ we define 
the determinant $\det(I+K)_{L^2(\Omega,\mu)}$ via a Fredholm expansion,
\begin{equation}\label{2.04}
\det(I+K)_{L^2(\Omega,\mu)}=\sum_{m=0}^\infty \frac 1{m!}\int_{\Omega^m}
\det (K(x_i,x_j))_{i,j=1,\dots,m} d^m\mu(x).
\end{equation}
We assume that $K$ is such that all the integrals are well-defined and
the series converges. For example, by Hadamard's inequality, it is sufficient
to require that $|K(x,y)|\le a(x)b(y)$, where $a\in L^1(\Omega,\mu)$,
$b\in L^\infty(\Omega,\mu)$. Note that if $\Omega=\{1,\dots,n\}$ and 
$\mu$ is counting measure this is the ordinary determinant 
$\det(\delta_{ij}+K(i,j))_{i,j=1,\dots,n}$. Let $K(x,y)$ be an integral kernel
from $L^2(\Omega_1,\mu_1)$ to $L^2(\Omega_2,\mu_2)$ and $L(x,y)$ an 
integral kernel from $L^2(\Omega_2,\mu_2)$ to $L^2(\Omega_1,\mu_1)$. Then
\begin{equation}
L\ast K(x,y)=\int_{\Omega_2} L(x,z)K(z,y)d\mu_2(z)
\notag
\end{equation}
is an integral kernel on $L^2(\Omega_1,\mu_1)$. Furthermore,
\begin{equation}\label{2.05}
\det(I+L\ast K)_{L^2(\Omega_1,\mu_1)}=\det(I+K\ast L)_{L^2(\Omega_2,\mu_2)}.
\end{equation}
This is easy to see using the Heine identity in the definition (\ref{2.04}).

Set
\begin{align}
b(i;u,\xi)&=\sum_{k=1}^n(A^{-1})_{ik}\phi_{-M,u}(x_k^{-M},\xi)\notag\\
c(u,\xi;j)&=\int_{\Lambda_m}\sum_{\ell=1}^{2M-1}\psi^{\ast(\ell-1)}(u,\xi;
v,\eta)g(v,\eta)\phi_{v,M}(\eta,x_j^M)d\nu(v,\eta),
\notag
\end{align}
so that, by (\ref{2.03}) and (\ref{2.05}),
\begin{align}
\frac{Z_{n,M}[g]}{Z_{n,M}[0]}&=
\det(\delta_{ij}+(b\ast c)(i,j))_{1\le i,j\le n}\notag\\
&=\det(I+c\ast b)_{L^2(\Lambda_M,\nu)}.
\notag
\end{align}
Now, a computation shows that
\begin{equation}
(c\ast b)(u,\xi;v,\eta)=(\sum_{\ell=1}^{2M-1}\psi^{\ast(\ell-1)})\ast
(g\tilde{K})(u,\xi;v,\eta),
\notag
\end{equation}
where $\tilde{K}$ is defined by (\ref{0.6}). Thus,
\begin{equation}\label{2.06}
\frac{Z_{n,M}[g]}{Z_{n,M}[0]}=\det( I+
(\sum_{\ell=1}^{2M-1}\psi^{\ast(\ell-1)})\ast
(g\tilde{K}))_{L^2(\Lambda_M,\nu)}.
\end{equation}
The kernel in (\ref{2.06}) has finite-rank so the sum (\ref{2.04}) in the
definition of the determinant actually has finitely many terms. 
We now claim that the right hand side of (\ref{2.06}) equals
\begin{equation}
\det(I-\psi+g\tilde{K})_{L^2(\Lambda_M,\nu)},
\notag
\end{equation}
which is what we want. Formally the computation goes as follows. 
The expression in (\ref{2.06}) is $\det(I-(I-\psi)^{-1}g\tilde{K})$
and we multiply this by $\det(I-\psi)=1$. Since we are only working 
with determinants defined by a Fredholm expansion the product rule
is not obvious, so we will give a proof in this special case.

Write $a=g\tilde{K}$. We will prove that for any $z,w\in\mathbb{C}$,
\begin{equation}\label{2.07}
\det(I+w\sum_{j=1}^mz^j\psi^{\ast(j-1)}\ast a)_{L^2(\Lambda_M)}=
\det(I-z\psi+zwa)_{L^2(\Lambda_M)},
\end{equation}
where $m=2M-1$. The left hand side is a polynomial in $z,w$ so it 
suffices to prove (\ref{2.07}) for $|z|,|w|$ sufficiently small. 
In that case, under our assumption on the $\phi_{r,r+1}$, all the 
expressions beloware well-defined and convergent. Write $b=
\sum_{j=1}^mz^j\psi^{\ast(j-1)}\ast a$. Then, see e.g. \cite{Me},
\begin{equation}
\det(I+wb)_{L^2(\Lambda_M)}=
\exp(\sum_{k=1}^\infty \frac{(-1)^{k+1}w^k}{k} \int_{\Lambda_M}
b^{\ast k}(t,t)d\nu(t))
\notag
\end{equation}
and
\begin{equation}
\det(I+z(-\psi+wa))_{L^2(\Lambda_M)}=
\exp(\sum_{k=1}^\infty \frac{(-1)^{k+1}z^k}{k} \int_{\Lambda_M}
(-\psi+wa)^{\ast k}(t,t)d\nu(t)).
\notag
\end{equation}
Set $c=za$, $d=z\psi$. It suffices to show that
\begin{align}\label{2.08}
&\sum_{k=1}^\infty \frac{(-1)^{k+1}w^k}{k} \int_{\Lambda_M}
\left(\sum_{j=0}^{m-1}d^j\ast c\right)^{\ast k} (t,t)d\nu(t)\notag\\
&=\sum_{k=1}^\infty \frac{(-1)^{k+1}}{k}\int_{\Lambda_M}
(-d+wc)^{\ast k}(t,t)d\nu(t).
\end{align}
The equality (\ref{2.08}) holds for $w=0$ since
$$
\int_{\Lambda_M}(-d)^{\ast k}(t,t)d\nu(t)=(-z)^k
\int_{\Lambda_M}\psi^{\ast k}(t,t)d\nu(t)=0
$$
if $k\ge 1$. This follows from $\phi_{r,s}=0$ for $r\ge s$. Hence it is
enough to show that the derivatives of the two sides of (\ref{2.08}) coincide,
\begin{align}
&\sum_{k=0}^\infty (-1)^kw^k\int_{\Lambda_M}
\left(\sum_{j=0}^{m-1}d^j\ast c\right)^{\ast (k+1)} (t,t)d\nu(t)\notag\\
&=\sum_{n=0}^\infty (-1)^n\int_{\Lambda_M}
((-d+wc)^{\ast n}\ast c)(t,t)d\nu(t).
\notag
\end{align}
To prove this last equality is a straightforward but somewhat tedious 
computation, which is based on expanding both sides and showing that the
coefficient of $w^k$ is the same on both sides. We omit the details.

We have proved
\begin{proposition}\label{P2.01}
Let $p_{n,M}(\bar{x})$ be defined by (\ref{0.3}) and assume that 
$\phi_{r,r+1}(x,y)$ satisfies $|\phi_{r,r+1}(x,y)|\le c(x)d(y)$ with $c
\in L^1(\mathbb{R},\mu)$, $d\in L^\infty(\mathbb{R},\mu)$, $-M\le r<M$. 
Furthermore, let $\Lambda_M=\{-M+1,\dots,M-1\}\times\mathbb{R}$, $\nu=
\lambda\otimes\mu$, where $\lambda$ is counting measure on 
$\{-M+1,\dots,M-1\}$ and let $g:\Lambda_M\to\mathbb{C}$ be a bounded function.
Then
\begin{equation}\label{2.09}
\int_{(\mathbb{R}^n)^{2M-1}}
\prod_{|\mu|<M}\prod_{j=1}^n(1+g(x_j^\mu))
p_{n,M}(\bar{x})d\mu(\bar{x})=\det(I+gK)_{L^2(\Lambda_M,\nu)},
\end{equation}
where $K$ is given by (\ref{0.5}), and the determinant is defined by using 
the Fredholm expansion (\ref{2.04}).
\end{proposition}

Theorem \ref{T0.1} is a direct consequence of (\ref{2.09}), 
compare with the discussion in \cite{TW2}.

If $X_m=\{1,\dots,m\}$ and $\lambda$ is counting measure on $X_m$, then
$L^2(X_m,\lambda)\cong\mathbb{R}^m$ and we have a chain of isomorphisms
$L^2(X_m\times\Omega,\lambda\otimes\mu)\cong
L^2(X_m,\lambda)\otimes L^2(\Omega,\mu)\cong
\mathbb{R^m}\otimes L^2(\Omega,\mu)\cong
L^2(\Omega,\mu)\oplus\dots\oplus L^2(\Omega,\mu)$, where we have $m$ terms 
in the last direct sum. We can think of an element in 
$\mathbb{R^m}\otimes L^2(\Omega,\mu)$ as a column vector $(f_1(x)\dots f_m(x)
)^t$, where $f_i\in L^2(\Omega,\mu)$, $1\le i\le m$. Hence, an operator on
$L^2(X_m\times\Omega,\lambda\otimes\mu)$ defined by an intgral kernel
$K(r,\xi;r',\xi')$ can be thought of as a block operator on these column
vectors with block kernel $(K(r,\xi;r',\xi'))_{1\le r,r'\le m}$.

We also want to prove Proposition \ref{P0.3}. Let us write $T_n(a)$
for the $n\times n$ Toeplitz matrix with symbol $a$ and $T(a)$ for the
one-sided infinite Toeplitz matrix with symbol $a$. Consider the
function $\tilde{K}^{n,M}(z,w)$ defined by (\ref{0.12}) and let the symbol
$a$ be given by (\ref{0.11}). Then,
\begin{align}\label{2.4}
\tilde{K}^{n,M}(z,w)&=\sum_{x,y\in\mathbb{Z}}\left( \sum_{i,j=1}^n
  \phi_{r,M}(x,1-i)[T_n^{-1}(a)]_{ij} \phi_{-M,s}(1-j,y)\right)
  z^xw^{-y} \notag \\
&=\sum_{i,j=1}^n\left(\sum_{x\in \mathbb{Z}}
  \hat{f}_{r,M}(1-i-x)z^{x+i-1}\right) z^{1-i}[T_n^{-1}(a)]_{ij} w^{j-1}
\notag\\  
&\times\left(\sum_{y\in \mathbb{Z}}
  \hat{f}_{-M,s}(y+j-1)w^{-y+1-j}\right) \notag \\
&= \frac zw f_{r,M}(\frac 1z)f_{-M,s}(\frac 1w)\sum_{i,j=1}^n z^{-i} 
[T_n^{-1}(a)]_{ij} w^j.
\end{align}
To proceed we need a formula for the inverse of a Toeplitz matrix. We
will use the following result which follows from theorem 1.15 and
theorem 2.15, together with its proof, in \cite{BoSi}.

\begin{proposition}\label{P2.4}
Assume that $a(z)=a_+(z)a_-(z)$, $z\in\mathbb{T}$, where
\begin{equation}\label{2.5}
a^+(z)=\sum_{n=0}^\infty a_n^+ z^n, \quad 
a^-(z)=\sum_{n=0}^\infty a_{-n}^- z^{-n},
\end{equation}
$\sum_{n=0}^\infty(|a_n^+|+|a_{-n}^-|)<\infty$, and that $(a(z)$ has winding
number zero. Furthermore, suppose that
\begin{equation}\label{2.6}
\sum_{n\in\mathbb{Z}} |n|^{\alpha}|\hat{a}_n|<\infty
\end{equation}
for some $\alpha>0$, where $\hat{a}_n$ is the Fourier coefficient of
$a(z)$. Using (\ref{2.5}) we can extend $a_+(z)$ to $|z|\le 1$ and
$a_-(z)$ to $\{|z|\ge 1\}\cup\{\infty\}$ and we assume that they have no
zeros in these regions. Then, $T_n(a)$ is invertible for $n$
sufficiently large and there is a constant $C$ (which depends on $a$)
such that 
\begin{equation}\label{2.7}
\left| [T_n^{-1}(a)]_{jk}-[T(a_+^{-1})T(a_-^{-1})]_{jk}\right|\le
C\min(\frac 1{(n+1-k)^\alpha}, \frac 1{(n+1-j)^\alpha})
\end{equation}
for $1\le j,k\le n$.
\end{proposition}

We can now give the proof of proposition \ref{P0.3}.

\begin{proof} ({\it  of Proposition \ref{P0.3}}). The
function $a$ defined by (\ref{0.11}) has a Wiener-Hopf factorization
$a=a_+a_-$ where
\begin{equation}
a_{\pm}(z)=\prod_{t=-M}^{M-1}f_t^{\pm}(z),
\notag
\end{equation}
and all the assumptions of the previous theorem are satisfied. By
(\ref{2.7}) 
\begin{align}
&\left|\sum_{i,j=1}^nz^{-i}[T_n^{-1}(a)]_{ij}w^j-
\sum_{i,j=1}^n z^{-i}[T(a_+)T(a_-^{-1})]_{ij}w^j\right|\notag\\
&\le C\sum_{i,j=1}^n |z|^{-i}|w|^{j}\min(\frac 1{(n+1-i)^\alpha},
\frac 1{(n+1-j)^\alpha})
\notag\\
&\le \frac C{(|z|-1)(1-|w|)} (\frac
1{n^\alpha}+|w|^{n/2}+ \frac 1{|z|^{n/2}}).\notag
\end{align}
Also,
\begin{equation}
\left|\sum_{\text{$i>n$ or
      $j>n$}}z^{-i}[T(a_+^{-1})T(a_-^{-1})]_{ij}w^j\right| \le
      C\frac{|w|^n+|1/z|^n}{(|z|-1)(1-|w|)}.
\notag
\end{equation}
Set $b_\pm=1/a_\pm$ and note that $(\hat{b}_+)_k=0$ if $k<0$ and 
$(\hat{b}_-)_k=0$ if $k>0$. We can now compute
\begin{align}
&\sum_{i,j=1}^\infty z^{-i}[T(a_+^{-1})T(a_-^{-1})]_{ij}w^j
 \notag \\     &=\sum_{k=1}^\infty\left(\sum_{i\in\mathbb{Z}} z^{-i+k}
(\hat{b}_+)_{i-k}\right)
\left(\sum_{j\in\mathbb{Z}} w^{j-k}
(\hat{b}_-)_{k-j}\right)\left(\frac wz\right)^k\notag\\
&=\frac 1{a_+(1/z)a_-(1/w)}\frac{w/z}{1-w/z}.
\notag
\end{align}
It follows that
\begin{align}
&\frac zwf_{r,M}(\frac 1z)f_{-M,s}(\frac 1w)
\sum_{i,j=1}^n z^{-i}[T(a_+^{-1})T(a_-^{-1})]_{ij}w^j\notag\\
&=\frac z{z-w}\frac{\prod_{t=r}^{M-1} f_t^+(\frac 1z)f_t^-(\frac 1z)
\prod_{t=-M}^{s-1}f_t^+(\frac 1w)f_t^-(\frac 1w)}{\prod_{t=-M}^{M-1}
f_t^+(\frac 1z)f_t^-(\frac 1w)}=\mathcal{K}^M_{r,s}(z,w),
\notag
\end{align}
and the proposition is proved.
\end{proof}

We have
\begin{equation}\label{2.10}
K(r,x;s,y)=\tilde{K}(r,x;s,y)=\frac 1{(2\pi i)^2}\int_{\gamma_{r_2}}
\frac{dz}{z^{x+1}}\int_{\gamma_{r_1}} w^{y-1}dw\frac z{z-w}G(z,w).
\end{equation}
if $r\ge s$, where $1-\epsilon<r_1<r_2<1+\epsilon$. Using the residue
theorem it follows that for $r<s$,
\begin{equation}\label{2.11}
K(r,x;s,y)=\frac 1{(2\pi i)^2}\int_{\gamma_{r_1}}
\frac{dz}{z^{x+1}}\int_{\gamma_{r_2}} w^{y-1}dw\frac z{z-w}G(z,w).
\end{equation}
$1-\epsilon<r_1<r_2<1+\epsilon$, compare \cite{OkRe}.

\subsection{The extended Airy kernel}
The extended Airy kernel is defined by (\ref{0.24}). We can also 
define a modification by 
\begin{equation}\label{2A.0}
\tilde{A}(\tau,\xi;\tau',\xi')= \int_0^\infty e^{-\lambda(\tau-\tau')}
\Ai(\xi+\lambda)\Ai(\xi'+\lambda)d\lambda,
\end{equation}
which is well-defined both for $\tau\ge \tau'$ and for $\tau<\tau'$ by
the following standard estimate for the Airy function,
\begin{equation}\label{2A.1}
|\Ai(\xi)|\le C_Me^{-2|\xi|^{3/2}/3}
\end{equation}
for $\xi\ge -M$. Both $A$ and $\tilde{A}$ have a useful double integral 
formula.

\begin{proposition}\label{L4.1}
The extended Airy kernel (\ref{0.24}) is given by
\begin{equation}\label{2A.2}
A(\tau,\xi;\tau',\xi')=-\frac 1{4\pi^2}\int_{\im z=\eta}dz
\int_{\im w=\eta'}dw\frac{e^{i\xi z+i\xi'
    w+i(z^3+w^3)/3}}{\tau'-\tau+i(z+w)}, 
\end{equation}
where $\eta,\eta'>0$ and $\eta+\eta'+\tau-\tau'<0$ in case
$\tau'>\tau$. Also, the modified kernel $\tilde{A}$, (\ref{2A.0}), is given
by the same formula but where we now require that
$\eta+\eta'+\tau-\tau'>0$.
\end{proposition}

\begin{proof}
This is straightforward using the identities
\begin{equation}
\int_0^\infty e^{-\lambda(\tau-\tau'-iz-iw)}d\lambda=-\frac
1{\tau'-\tau+i(z+w)} 
\notag
\end{equation}
if $\tau-\tau'+\eta+\eta'>0$ and
\begin{equation}
\int_0^\infty e^{-\lambda(\tau'-\tau+iz+iw)}d\lambda=\frac
1{\tau'-\tau+i(z+w)} 
\notag
\end{equation}
if $\tau-\tau'+\eta+\eta'<0$.
\end{proof}

If we move the contour of integration between the two cases in proposition
\ref{L4.1} we pick up a contribution from the singularity and we obtain
\begin{equation}\label{2A.3}
A(\tau,\xi;\tau',\xi')=\tilde{A}(\tau,\xi;\tau',\xi')-\phi_{\tau,\tau'}
(\xi,\xi'),
\end{equation}
where $\phi_{\tau,\tau'}\equiv 0$ if $\tau\ge\tau'$ and
\begin{equation}
\phi_{\tau,\tau'}(\xi,\xi')=\frac 1{\sqrt{4\pi (\tau'-\tau)}}
  e^{-(\xi-\xi')^2/4(\tau'-\tau)
  -(\tau'-\tau)(\xi+\xi')/2+(\tau'-\tau)^3/12}
\end{equation}
if $\tau<\tau'$. Combining (\ref{0.24}), (\ref{2A.0}) and (\ref{2A.3}) 
we see that
\begin{equation}\label{2A.4}
\phi_{\tau,\tau'}(\xi,\xi')=
\int_
{-\infty}^\infty e^{-\lambda(\tau-\tau')}
\Ai(\xi+\lambda)\Ai(\xi'+\lambda)d\lambda
\end{equation}
if $\tau<\tau'$. We would also like to show that the operator in 
(\ref{0.28}) is actually a trace class operator.

\begin{proposition}\label{PA.1}
Let $f(\tau,x)$ bew a non-negative function in $L^\infty(\mathbb{R})$ for
each $\tau\in\{\tau_1,\dots,\tau_m\}$, where $\tau_1<\dots<\tau_m$. Assume
also that $f(\tau_k,x)=0$ if $x<M_k$ for some number $M_k$, $k=1,\dots,m$.
Then, the kernel
\begin{equation}
f(\tau,x)^{1/2}A(\tau,x;\tau',x')f(\tau',x')^{1/2}
\notag
\end{equation}
defines a trace class operator on $L^2(\{\tau_1,\dots,\tau_m\}\times
\mathbb{R})$, where we have counting measure $\lambda$ on $
\{\tau_1,\dots,\tau_m\}$ and Lebesgue measure $\mu$ on $\mathbb{R}$.
\end{proposition}

\begin{proof}
We will prove the result by factoring into two Hilbert-Schmidt operators.
Let $H(t)=1$ if $t<0$ and $H(t)=0$ if $t\ge 0$. Set
\begin{equation}
\tilde{B}(\tau,x;\tau',x')=H(\tau-\tau')\int_{-\infty}^\infty
e^{-y(\tau-\tau')}\Ai(x+y)\Ai(x'+y)dy.
\notag
\end{equation}
For $i<j$ we define
\begin{equation}
\tilde{B}_{ij}(\tau,x;\tau',x')=\tilde{B}(\tau,x;\tau',x')\delta_{\tau,\tau_i}
\delta_{\tau',\tau_j}
\notag
\end{equation}
so that
\begin{equation}
\tilde{B}(\tau,x;\tau',x')=\sum_{1\le i<j\le m}
\tilde{B}_{ij}(\tau,x;\tau',x').
\notag
\end{equation}
Since, by (\ref{2A.3}) and (\ref{2A.4}), $A=\tilde{A}-\tilde{B}$, it suffices
to show that $f^{1/2}\tilde{A}f^{1/2}$ and $f^{1/2}\tilde{B_{ij}}f^{1/2}$,
$1\le i<j\le m$, are trace class operators.

Set
\begin{equation}
a(\tau,x;\sigma,y)=\frac 1{\sqrt{m}}f(\tau,x)^{1/2}\Ai(x+y)e^{-y(\tau-\sigma)}
\chi_{[0,\infty)}(y)
\notag
\end{equation}
\begin{equation}
b(\sigma,y;\tau',x')=\frac 1{\sqrt{m}}\chi_{[0,\infty)}(y)
\Ai(x'+y)e^{-y(\sigma-\tau')}f(\tau',x')^{1/2}
\notag
\end{equation}
Then $a $ and $b$ are Hilbert-Schmidt kernels on $L^2(\Lambda_m,\lambda
\otimes\mu)$, $\lambda_m=\{\tau_1,\dots,\tau_m\}\times\mathbb{R}$. We have
\begin{align}
&\int_{\Lambda_m}\int_{\Lambda_m}|a(\tau,x;\sigma,y)|^2d(\lambda
\otimes\mu)(\tau,x)d(\lambda\otimes\mu)(\sigma,y)\notag\\
&=\frac 1m\int_{\Lambda_m}\int_{\Lambda_m}f(\tau,x)\Ai(x+y)^2
\chi_{[0,\infty)}(y)e^{-2y(\tau-\sigma)}dxdyd\lambda(\tau)d\lambda(\sigma)
\notag\\
&\le\frac{||f||_\infty}m \sum_{i,j=1}^m\int_{M}^\infty dx\int_0^\infty dy
\Ai(x+y)^2 e^{2y(\tau_m-\tau_1)},
\notag
\end{align}
where $M=\min(M_1,\dots,M_m)$. Using (\ref{2A.1}) we see that the integral 
in the last expression is $<\infty$. The proof that $b$ is a
 Hilbert-Schmidt kernel is analogous. Now,
\begin{align}
&\int_{\Lambda_m}a(\tau,x;\sigma,y) b(\sigma,y;\tau',x')d(\lambda\otimes\mu)
(\tau,y)\notag\\
&\frac 1m\sum_{\sigma\in\{\tau_1,\dots,\tau_m\}} f(\tau,x)^{1/2}
f(\tau',x')^{1/2}\int_0^\infty e^{-y(\tau-\tau')}\Ai(x+y)\Ai(x'+y)dy\notag\\
&=f(\tau,x)^{1/2}\tilde{A}(\tau,x;\tau',x')f(\tau',x')^{1/2}.
\notag
\end{align}
Hence, the operator $f^{1/2}\tilde{A}f^{1/2}$ is trace class.

Next, set
\begin{equation}
c_{ij}(\tau,x;\tau',x')=\frac 1{\sqrt{m}} f(\tau,x)^{1/2}\Ai(x+y) 
e^{-y(\tau-\tau_j)/2}\delta_{\tau,\tau_i}
\notag
\end{equation}
(it is independent of $\sigma$) and
\begin{equation}
d_{ij}(\tau,x;\tau',x')=\frac 1{\sqrt{m}} f(\tau',x')^{1/2}\Ai(x'+y) 
e^{-y(\tau_i-\tau')/2}\delta_{\tau',\tau_j}.
\notag
\end{equation}
Then,
\begin{align}
&\int_{\Lambda_m} c_{ij}(\tau,x;\tau',x')d_{ij}(\tau,x;\tau',x')
d(\lambda\otimes\mu)(\sigma,y)\notag\\
&=\frac 1m\sum_{\sigma\in\{\tau_1,\dots,\tau_m\}} f(\tau,x)^{1/2}
f(\tau',x')^{1/2}\int_{-\infty}^\infty e^{-y(\tau-\tau')}\Ai(x+y)\Ai(x'+y)dy
\delta_{\tau,\tau_i}\delta_{\tau',\tau_j}\notag\\
&=\tilde{B}_{ij}(\tau,x;\tau',x').
\notag
\end{align}
It remains to prove that $c_{ij}$ and $d_{ij}$ are Hilbert-Schmidt kernels. 
Consider $c_{ij}$; the proof for $d_{ij}$ is similar. We get
\begin{align}
&\frac 1m\sum_{\sigma,\tau\in\{\tau_1,\dots,\tau_m\}} 
\int_{-\infty}^\infty\int_{-\infty}^\infty f(\tau,x)\Ai(x+y)^2
e^{-y(\tau-\tau_j)}dxdy\delta_{\tau,\tau_i}\notag\\
&=\int_{-\infty}^\infty\int_{-\infty}^\infty f(\tau_i,x)\Ai(x+y)^2
e^{-y(\tau_i-\tau_j)}dxdy\notag\\
&\le ||f||_\infty\int_{M_i}dx\int_{-\infty}^\infty dy\Ai(x+y)^2
e^{y(\tau_j-\tau_i)}\notag\\
&\le ||f||_\infty\int_{M_i}dx\int_{0}^\infty dy\Ai(x+y)^2
e^{y(\tau_j-\tau_i)}\notag\\
&+||f||_\infty\int_{M_i}dx\int_{-\infty}^0 dy\Ai(x+y)^2
e^{y(\tau_j-\tau_i)}.
\notag
\end{align}
The first integral in the last expresion
is $<\infty$ by (\ref{2A.1}). Now, by (\ref{2A.1})
\begin{equation}
\int_{M_i}^\infty \Ai(x+y)^2dx=\int_{M_i+y}^\infty \Ai(x)^2dx\le C(1+|y|)
\notag
\end{equation}
since the Airy function is bounded. Hence, 
\begin{align}
&||f||_\infty\int_{M_i}dx\int_{-\infty}^0 dy\Ai(x+y)^2
e^{y(\tau_j-\tau_i)} \notag\\
&\le C||f||_\infty \int_{-\infty}^0(1+|y|)e^{y(\tau_j-\tau_i)}dy<\infty,
\notag
\end{align}
since $\tau_j-\tau_i>0$. This completes the proof.
\end{proof}

\subsection{An example: random walks on the discrete circle}
We will consider non-intersecting walks on the set $\mathbb{Z}_N$ of
integers modulo $N$, the discrete circle. 
This type of model has been analyzed in \cite{For} and we will show
how it fits into the present formalism.
We have $2M-1$ copies of 
$\mathbb{Z}_N$, where the first and the last are identified so that we
have periodic boundary conditions in the time direction. We will
have are non-=intersecting paths on the discrete torus. Let
$x^r\in\mathbb{Z}_N^n$ be the particle configuration ($n$ particles)
on the $r$:th discrete circle, $|r|<M$,$x^{-M+1}\equiv
x^{M-1}$. Assume that $n$ is odd, $n=2\nu+1$ and that the
transition probabilities for the walks are given by
\begin{equation}
\phi_{r,r+1}(x,y)=\begin{cases}
p, &\text{if $y-x=1$}\\
q,  &\text{if $y-x=0$.}\\
0. &\text{otherwise}
\end{cases}
\notag
\end{equation}
for $x,y\in\mathbb{Z}_N$, where $p,q\ge 0$ and $p+q=1$. The transition
probability for non-intersecting paths from a configuration $x^r$ to a
configuration $x^{r+1}$ is
\begin{equation}
\det(\phi_{r,r+1}(x_i^r,x_j^{r+1})_{1\le i,j\le n}).
\notag
\end{equation}
Write $\bar{x}=(x^{-M+1},\dots,x^{M-1})\in(\mathbb{Z}_N)^{2M-1}$ for
the total configuration. The probability of $\bar{x}$ is
\begin{equation}\label{2.101}
q_{n,N,M}(\bar{x})=\prod_{r=-M+1}^{M-2} 
\det(\phi_{r,r+1}(x_i^r,x_j^{r+1})_{1\le i,j\le n}.
\end{equation}
We will use discrete Fourier series on $\mathbb{Z}_N$,
\begin{equation}
\hat{f}(n)=\frac 1N \sum_{\ell\in\mathbb{Z}_N} f(\ell) z^{-\ell
  n},\quad
f(\ell)=\sum_{n\in\mathbb{Z}_N}\hat{f}(n) z^{\ell n},
\notag
\end{equation}
where $z=e^{2\pi i/N}$. Also, we can represent Kronecker's delta on
$\mathbb{Z}_N$ as
\begin{equation}\label{2.102}
\delta_{xy}=
\frac 1N \sum_{k\in\mathbb{Z}_N} z^{k(x-y)}.
\end{equation}
Let $\tilde{\mathbb{Z}}_N^n=\{x\in\mathbb{Z}_N^n\,;\, 0\le
x_1<\dots<x_n<N\}$ be all ordered configurations of $n$ particles on
$\mathbb{Z}_N$ . If $x^{-M+1}, x^{M-1}\in\tilde{\mathbb{Z}}_N^n$, then
\begin{equation}
\det (\delta_{x_i^{-M+1},x_j^{M-1}})_{1\le i,j\le n} 
=\delta_{x^{-M+1},x^{M-1}}
\notag
\end{equation}
by the definition of the determinant. This determinant can be
rewritten using (\ref{2.102}) and Heine's identity,
\begin{align}\label{2.103}
&\delta_{x^{-M+1},x^{M-1}}=\det(\frac 1N\sum_{k\in\mathbb{Z}_N}
z^{k(x_i^{-M+1} -x_j^{M-1})})_{1\le i,j\le n}\\
&=\frac 1{n!N^n}\sum_{k_1,\dots,k_n\in\mathbb{Z}_n} 
\det(z^{k_ix_j^{-M+1}})_{1\le i,j\le n}
\det(z^{-k_ix_j^{M-1}})_{1\le i,j\le n},
\notag
\end{align}
if $x^{-M+1},x^{M-1}\in\tilde{\mathbb{Z}}_N^n$. This leads us to the measure
\begin{equation}
p_{n,N,M}(\bar{x})=\frac 1{(n!)^{2M-1}Z_{n,N,M}}
q_{n,N,M}(\bar{x})\delta_{x^{-M+1},x^{M-1}} 
\notag
\end{equation}
where $Z_{n,N,M}$ is the normalization constant. 

Let $g(r,x)$, $|r|<M$, $x\in\mathbb{Z}_N$, be given functions
and set
\begin{equation}
G(\bar{x})=\prod_{|r|<M}\prod_{j=1}^n (1+g(r,x_j^r)).
\notag
\end{equation}
We want to compute the expectation 
\begin{align}\label{2.104}
&\sum_{\bar{x}\in(\mathbb{Z}_N^n)^{2M-1}} G(\bar{x}
)p_{n,N,M}(\bar{x})\\
&=\frac{n!}{N^nZ_{n,N,M}}\sum_{\bar{x}\in(\tilde{\mathbb{Z}}_N^n)^{2M-1}}
\sum_{k\in\tilde{\mathbb{Z}}_N^n} G(\bar{x}) w_{n,N,M}(k;\bar{x}),
\notag
\end{align}
where
\begin{align}\label{2.105}
w_{n,N,M}(k;\bar{x})&=\det(z^{k_ix_j^{-M+1}}) q_{n,N,M}(\bar{x})
\det(z^{-k_ix_j^{M-1}})\\
&=\prod_{r=-M}^{M-1}\det(\phi_{r,r+1}(x_i^r,x_j^{r+1})).
\notag
\end{align}
Here we have set $\phi_{-M,-M+1}(k_i,x_j^{-M+1})=z^{k_ix_j^{-M+1}}$,
$\phi_{M-1,M}(x_i^{M-1},k_j)=z^{-k_jx_i^{M-1}}$ and
$x_i^{-M}=x_i^M=k_i\in\mathbb{Z}_N$. We have a measure of the form
(\ref{0.1}).
Set
\begin{equation}\label{2.106}
Z_{n,N,M}(k)=\sum_{\bar{x}\in(\tilde{\mathbb{Z}}_N^n)^{2M-1}} 
w_{n,N,M}(k;\bar{x})
\end{equation}
and note that $G\equiv 1$ in (\ref{2.104}) gives
\begin{equation}\label{2.107}
Z_{n,N,M}=\frac {n!}{N^n}\sum_{k\in\tilde{\mathbb{Z}}_N^n} Z_{n,N,M}(k).
\end{equation}
Let us also write
\begin{equation}
p_{n,N,M}(k,\bar{x})=\frac
1{(n!)^{2M-1}Z_{n,N,M}(k)}w_{n,N,M}(k;\bar{x}).
\notag
\end{equation}
The expectation (\ref{2.104}) can then be written, using
(\ref{2.107}),
\begin{equation}\label{2.108}
\sum_{k\in\tilde{\mathbb{Z}}_N^n} \frac{Z_{n,N,M}(k)}{
\sum_{k\in\tilde{\mathbb{Z}}_N^n}Z_{n,N,M}(k)} E_{n,N,M}(k;G)
\end{equation}
where $E_{n,N,M}(k;G)$ is the ``expectation'' of $G$ with repect to the
measure $p_{n,N,M}(k,\bar{x})$. This ``expectation'' can be computed
using the standard framework. Let $\hat{f}(n)$ be equal to $q$ if
$n=0$, $p$ if $n=1$ and 0 otherwise, $n\in\mathbb{Z}_N$, so that
$\phi_{r,r+1}(x,y)=\hat{f}(y-x)$, for $-M<r<M-1$. Then
$f(\ell)=q+pz^\ell$, $\ell\in\mathbb{Z}_N$, and by standard properties
of convolution
\begin{equation}\label{2.109}
\phi_{r,s}(x,y)=\hat{f}^{s-r}(y-x),
\end{equation}
$-M<r<s<M-1$. From this we see that 
\begin{align}
A_{ij}&=\phi_{-M,M}(k_i,k_j)=\sum_{x,y\in\mathbb{Z}_N} z^{k_ix}
\hat{f}^{2M-2} (y-x) z^{-k_j y}\notag\\
&=Nf(N-k_i)^{2M-2}\delta_{k_i,k_j}=Nf(N-k_i)^{2M-2}\delta_{ij}
\notag
\end{align}
if $k\in \tilde{\mathbb{Z}}_N^n$. Thus,
\begin{equation}\label{2.110}
A=(Nf(N-k_i)^{2M-2}\delta_{ij})_{i,j=1,\dots,n}.
\end{equation}
Now,
\begin{equation}\label{2.110'}
Z_{n,N,M}(k)=\det A=N^n\prod_{i=1}^n(q+pe^{-2\pi ik_i/N})^{2M-2}.
\end{equation}
This is always non-zero if $p\neq q$. If $p=q=1/2$, then we assume
that $N$ is odd, which also implies that $\det A\neq 0$. 
We obtain
\begin{align}\label{2.111}
\tilde{K}_n(k;r,x;s,y)&=\sum_{i,j=1}^n\left( \sum_{\ell\in
    \mathbb{Z}_N} \hat{f}^{M-r-1}(\ell-x)z^{-k_i\ell}\right)\frac 1N
    f(N-k_i)^{2-2M}\delta_{ij}
\\
&\left( \sum_{m\in
    \mathbb{Z}_N} \hat{f}^{s+M-1}(y-m)z^{k_jm}\right) \notag\\
&=\frac 1N\sum_{i=1}^nf(N-k_i)^{s-r}z^{k_i(y-x)},
\notag
\end{align}
where we have indicated the dependence of the kernel on $k$.
Note that this kernel is independent of $M$. We have
\begin{equation}\label{2.112}
K_n(k;r,x;s,y)=\frac 1N\sum_{i=1}^n f(N-k_i)^{s-r}z^{k_i(y-x)}
-\phi_{r,s}(x,y), 
\end{equation}
where
\begin{equation}\label{2.113}
\phi_{r,s}(x,y)=\frac 1N\sum_{\ell\in\mathbb{Z}_N}
f(\ell)^{s-r}z^{\ell (y-x)}
\end{equation}
if $s>r$ and $\phi_{r,s}(x,y)=0$ if $s\le r$.
Computations similar to those leading up to formula (\ref{3.32}) below
show that if we assume that $g(r,\cdot)\equiv 0$
if $|r|>M_0$, then 
$E_{n,N,M}(k,G)=E_{n,N,M_0}(k,G)$ for $M\ge M_0$. Hence, the
expectation (\ref{2.108}) can be written 
\begin{equation}\label{2.115}
\sum_{k\in\tilde{\mathbb{Z}}_N^n}\left(\lim_{M\to\infty}\frac
  {Z_{n,N,M}(k)}
  {\sum_{k\in\tilde{\mathbb{Z}}_N^n}Z_{n,N,M}(k)}\right) 
  E_{n,N,M_0}(k,G).
\end{equation}

\begin{lemma}\label{C1}
Let $\alpha_i=i-1$, $i=1,\dots,\nu+1$, $\alpha_{2\nu+2-i}=N-i$,
$i=1,\dots,\nu$, $n=2\nu+1$. If $k\in\tilde{\mathbb{Z}}_N^n$, then
\begin{equation}\label{2.116}
\lim_{M\to\infty}\frac
  {Z_{n,N,M}(k)}
  {\sum_{k\in\tilde{\mathbb{Z}}_N^n}Z_{n,N,M}(k)}
=\delta_{\alpha,k}.
\end{equation}
\end{lemma}
\begin{proof}
We use the explicit formula (\ref{2.110'}) for $Z_{n,N,M}(k)$,
\begin{equation}
Z_{n,N,M}(k)=N^n\prod_{j=1}^n(q+pe^{-2\pi i k_j/N})^{2M-2}.
\notag
\end{equation}
Now,
\begin{equation}
|q+pe^{-2\pi i k_i/N}|^2=p^2+q^2+2pq\cos\frac {2\pi ik_i}N.
\notag
\end{equation}
This is maximal $(=1)$ if $k_i=0 (=N)$ and it is easy to see that
\begin{equation}
\prod_{j=1}^n(p^2+q^2+2pq\cos\frac {2\pi ik_j}N)\le
\prod_{j=1}^n(p^2+q^2+2pq\cos\frac {2\pi i\alpha_j}N)
\notag
\end{equation}
with strict inequality unless $k=\alpha$. This completes the proof.
\end{proof}

Hence, if $g(r,\cdot)\equiv 0$ for $|r|\ge M_0$, then
\begin{equation}
\lim_{M\to\infty}\sum_{\bar{x}\in (\mathbb{Z}_N^n)^{2M-1}} G(\bar{x})
p_{n,N,M}(\bar{x})= E_{n,N,M_0}(\alpha;G).
\notag
\end{equation}
From this it follows that the correlation kernel $K(r,x;s,y)$  
on the cylinder $\mathbb{Z}\times \mathbb{Z}_N$ is given by
$K(\alpha;r,x;s,y)$. We obtain the following proposition.
\begin{proposition}\label{C2}
The correlation function for $n=2\nu+1$ non-intersecting walks on the
infinite cylinder $\mathbb{Z}\times \mathbb{Z}_N$ as defined above is
given by
\begin{align}\label{2.117}
K(r,x;s,y)&=\frac 1N\sum_{j=-\nu}^\nu (q+pe^{2\pi i j/N})^{s-r}
e^{2\pi i j(x-y)/N}\\
&-\omega_{r,s}\frac 1N\sum_{j=-\nu}^{N-\nu-1}(q+pe^{2\pi i j/N})^{s-r}
e^{2\pi i j(x-y)/N}
\notag
\end{align}
where $\omega_{r,s}=1$ if $r<s$, $\omega_{r,s}=0$ if $r\ge s$.
\end{proposition}

The induced measure on $\mathbb{Z}_N$ is given by 
\begin{align}
\frac 1{n!}\det (K(0,x_\mu;0,x_\nu))_{1\le\mu,\nu\le n}&= 
\frac 1{n!}\det (\frac 1N\sum_{j=-\nu}^{\nu} e^{2\pi i
  j(x_\mu-x_\nu)/N} )_{1\le\mu,\nu\le n}\notag\\
=\frac 1{n!N^n}\prod_{1\le\mu<\nu\le n} |e^{2\pi i x_\mu/N}-
e^{2\pi i x_\nu/N}|^2,
\notag
\end{align}
the equilibrium measure on $\mathbb{Z}_N$ (discrete CUE), see \cite{KOCR}.

We can take the limit $n,N\to\infty$, $n/N\to\rho$, $0<\rho<1$, and
obtain a limiting determinantal process on $\mathbb{Z}^2$.

\begin{proposition}\label{C3}
The correlation function for the determinantal process on
$\mathbb{Z}^2$ induced by non-intersecting random walks as defined
above is given by
\begin{equation}\label{1.119a}
K(r,x;s,y)=\int_{-\rho/2}^{\rho/2} (q+pe^{2\pi i\theta})^{s-r}
e^{2\pi i \theta (x-y)} d\theta
\end{equation}
if $r\ge s$, and 
\begin{equation}\label{1.119b}
K(r,x;s,y)=-\int_{\rho/2}^{1-\rho/2} (q+pe^{2\pi i\theta})^{s-r}
e^{2\pi i \theta (x-y)} d\theta
\end{equation}
if $r< s$ for $0<\rho<1$.
\end{proposition}

This kernel is related to the $B^\pm$-kernels in \cite{OkRe}. Compare
also with \cite{YNS}.

\section{Multi-layer discrete PNG}

We will discuss how the PNG model defined by (\ref{0.16}), in the case
when $\omega(x,t)$, $(x,t)\in\mathbb{Z}\times\mathbb{N}$, satisfies
$\omega(x,t)=0$ if $t-x$ is even or if $|x|>t$, can be embedded as
the top curve in a multi-layer process given by a family of 
non-intersecting paths. We think of the $\omega(x,t)$:s as given numbers.
The initial condition is $h(x,0)=0$, $x\in\mathbb{Z}$. We extend $h(x,t)$
to all $x\in\mathbb{R}$ by letting $h(x,t)=h([x],t)$, which makes 
it right continuous at the jumps. Note that it follows immediately that
$h(x,t)=0$ if $x<-t+1$ or $x>t$. 

For $t-x$ {\it odd} we define the {\it jumps},
\begin{align}\label{3.02}
\eta^+(x,t)&=h(x,t)-h(x-1,t)\\
\eta^-(x,t)&=h(x,t)-h(x+1,t).
\notag
\end{align}
We will see below that $\eta^+,\eta^-\ge0$ and we should think of 
$\eta^+(x,t)$ as a positive jump at $x$ at time $t$, and 
$\eta^-(x,t)$ as a the size of a negative jump at $x+1$ at time $t$.
Define
\begin{equation}\label{3.03}
T\omega(x,t)=\min(\eta^+(x+1,t-1),\eta^-(x-1,t-1))
\end{equation}
if $t-x$ is odd and $T\omega(x,t)=0$ if $t-x$ is even.
\begin{claim}\label{Cl3.1}
The jumps $\eta+$ and $\eta^-$ satisfy the following evolution equation
\begin{align}\label{3.04}
\eta^+(x+1,t+1)&=\max(\eta^+(x+2,t)-\eta^-(x,t),0)+\omega(x+1,t-1)\\
\eta^+(x+1,t+1)&=\max(\eta^-(x,t)-\eta^+(x+2,t),0)+\omega(x+1,t-1)
\notag
\end{align}
for $t-x$ odd. Furthermore $\eta^+(x,t)$ and $\eta^-(x,t)$ are $\ge 0$.
\end{claim}
\begin{proof}
We proceed by induction on $t$. Assume that $\eta^+(x,t),\eta^-(x,t)\ge 0$
for all $x$ such that $t-x$ is odd. We will prove that then (\ref{3.04})
holds, and hence $\eta^+(x+1,t+1),\eta^-(x+1,t+1)\ge 0$ 
for all $x$ such that $t-x$ is odd. Obviously our induction assumption
is true for $t=0$. Note that $h(x+1,t)=h(x,t)-\eta^-(x,t)$,
$h(x+2,t)=h(x,t)+\eta^+(x+2,t)-\eta^-(x,t)$ and $h(x-1,t)=
h(x,t)-\eta^+(x,t)$. It follows from (\ref{0.16}), our induction aasumption
and $\omega(x,t+1)=0$, that
\begin{align}
h(x+1,t+1)&=h(x,t)+\max(0,\eta^+(x+2,t)-\eta^-(x,t))+\omega(x+1,t+1)
\notag\\
h(x,t+1)&=h(x,t)
\notag
\end{align}
and the first half of (\ref{3.04}) follows. The proof of the second half
is analogous.
\end{proof}

There is also an inverse recursion formula.
\begin{claim}\label{Cl3.02}
If $t-x$ is odd, then
\begin{align}\label{3.05}
\omega(x+1,t+1)&=\min(\eta^-(x_+1,t+1),\eta^+(x_+1,t+1))\notag\\
\eta^+(x,t)&=\eta^+(x-1,t+1)-\omega(x-1,t+1)+T\omega(x-1,t+1)\\
\eta^-(x,t)&=\eta^+(x+1,t+1)-\omega(x+1,t+1)+T\omega(x+1,t+1)
\notag
\end{align}
\end{claim}
\begin{proof}
The first equation folows immediately from (\ref{3.04}). 
From (\ref{3.03}) and (\ref{3.04})  we see that the right hand side of the
second equation in (\ref{3.05}) equals
\begin{equation}
\max(\eta^+(x,t)-\eta^-(x-2,t),0)+\min(\eta^+(x,t),\eta^-(x-2,t)),
\notag
\end{equation}
which equals $\eta^+(x,t)$. The proof of the last equation is similar.
\end{proof}

From this claim we immediately deduce the following
\begin{claim}\label{Cl3.03}
If we know $\eta^+(x+1,t+1)$, $\eta^-(x+1,t+1)$ for all $x$ such that
$t-x$ is odd, and $T\omega(x,s)$ for $s\le t+1$ and all $x$, we can 
reconstruct $\omega(x,s)$, $s\le t+1$, $x\in\mathbb{Z}$, uniquely.
\end{claim}

Let a coordinate system $(i,j)$ be related to the $(x,t)$ coordinate
system via the transformation
\begin{equation}\label{3.06}
(x,t)=(i-j,i+j-1),
\end{equation}
and define $w(i,j)$ by (\ref{0.17}). Then $w(i,j)=0$ if $(i,j)\notin
\mathbb{Z}_+^2$, and this condition corresponds exactly to our assumptions
on $\omega(x,t)$. Similarly to (\ref{0.17}) we define
\begin{equation}\label{3.07}
Tw(i,j)=T\omega(i-j,i+j-1).
\end{equation}

\begin{claim}\label{Cl3.04}
Assume that $w(i,j)=0$ if $i$ or $j$ is $\le s$. Then, $Tw(i,j)=0$ if
$i$ or $j$ is $\le s+1$.
\end{claim}
\begin{proof}
It follows from the condition on $w(i,j)$ that $\omega(x,t)=0$ if 
$(t+x+1)/2\le s$ or $(t-x+1)/2\le s$, which implies, using (\ref{0.16}),
that $h(x,t)=0$ under the same conditions. It follows from (\ref{3.02}),
(\ref{3.03}) and (\ref{3.07}) that $Tw(i,j)=0$ if $h(i-j+1,i+j-2)=0$ or
$h(i-j-1,i+j-2)=0$. Now, $h(i-j+1,i+j-2)=0$ if $i\le s$ or $j\le s+1$, and
$h(i-j-1,i+j-2)=0$ if $i\le s+1$ or $j\le s$. Hence $Tw(i,j)=0$ if
$i\le s+1$ or $j\le s+1$.
\end{proof}

It follows from claim \ref{Cl3.04} that $T^nw(i,j)=0$ if $i$ 
or $j$ is $\le n$, since $w(i,j)=0$ if $i$ or $j$ is $\le 0$. 
Hence $T^n(\omega(x,t)-0$ if $t\le 2n-1$, since $i+j-1\le 2n-1$ 
implies that $i$ or $j$ is $\le n$. We formulate this as our next claim.

\begin{claim}\label{Cl3.05}
If $t\le 2n-1$, then $T^n\omega(x,t)=0$.
\end{claim}

Let $h_i(x,t)$, $i\ge 0$, be the PNG process defined by (\ref{0.16}) with
$\omega(x,t+1)$ replaced by $T^i\omega(x,t+1)$, and with initial condition 
$h_i(x,o)=-i$. We let $T^0\omega(x,t+1)=\omega(x,t+1)$, so $h_0(x,t)=
h(x,t)$ is our original growth process. It follows from claim \ref{Cl3.05} 
that at time $t=2n-1$ only $h_0,\dots, h_{n-1}$ can be non-trivial, i.e. 
$h_i(x,2n-1)=-i$ for all $x$ if $i\ge n$. Combining claim \ref{Cl3.03}
and claim \ref{Cl3.05} we get
\begin{claim}\label{Cl3.06}
Given $h_i(x, 2n-1)$, $x\in\mathbb{Z}$, $i=0,\dots,n-1$, we can uniquely
reconstruct $\{\omega(x,t)\,;\, t\le 2n-1,x\in\mathbb{Z}\}$.
\end{claim}

We can think of $h_i$ at time $2n-1$ as a directed path from $(-2n+1,-i)$ to
$(2n-1,-i)$ which has up-steps $\eta^+(2m,2n-1)$ at even $x$-coordinates,
$x=2m$ and down-steps $\eta^-(2m,2n-1)$ at odd $x$-coordinates, $x=2m+1$,
$|m|<n$, and horizontal steps in between. According to \ref{Cl3.06} there is a
bijection between these paths $h_0,\dots,h_{n-1}$ and the set $\{\omega (x,t)
\,;\,t\le 2n-1\, ,\,x\in\mathbb{Z}\}$. 
We set $h_i(x,t)=h_i([x],t)$ for $x\in\mathbb{R}$.
The paths obtained are nonintersecting:
\begin{claim}\label{Cl3.07}
If $t-x$ is odd, then
\begin{equation}\label{3.08}
h_{i+1}(x,t)<h_i(x-0,t)
\end{equation}
and if $t-x$ is even
\begin{equation}\label{3.09}
h_{i+1}(x-0,t)<h_i(x,t),
\end{equation}
so that corners will not meet.
\end{claim}
\begin{proof}
This is proved by induction on $t$. It is clearly true for $t-0$. 
If it is true at time $t$ it is still true after forming the maximum in
(\ref{0.16}) (deterministic step). (Note that $h_i$ and $h_{i+1}$ have
up-steps/down-steps at the same positions.) From the definition (\ref{3.03})
it is still true after adding $T^{i+1}\omega(x,t)$ to the lower curve.
\end{proof}

In order to understand how a geometric distribution (\ref{0.18}) on the 
$w(i,j)$ is transported to a measure on the non-imtersecting paths, 
we will assign weights to the jumps. Let $a_i$ and $b_j$ be given variables.
The jumps are assigned weights as follows: $\eta^+(x,t)$ has weight 
$a_i^{\eta^+(x,t)}$, $i=(t+x+1)/2$ and $\eta^-(x,t)$ has weight 
$b_j^{\eta^-(x,t)}$, $j=(t-x+1)/2$. Also, to $T^k\omega(x,t)$ we assign the
weight $(a_i,b_j)^{T^k\omega(x,t)}$ with the same correspondence between
$(i,j)$ and $(x,t)$, $k\ge 0$. The proof of the next claim is a 
straightforward computation using the defintions of the quantities involved.
\begin{claim}\label{Cl3.08}
The product of the weights of $\eta^+(x-1,t+1)$, $\eta^-(x-1,t+1)$ and
$T\omega(x-1,t+1)$ equals the product of the weights of $\eta^-(x-2,t)$,
$\eta^+(x,t)$ and $\omega(x-1,t+1)$.
\end{claim}

Using this claim we can show that the measure is transported in the way we 
want.

\begin{claim}\label{Cl3.09}
The product of all the weights of all the jumps in the multi-layer PNG,
$h_0,\dots,h_{n-1}$, at time $t=2n-1$ equals,
\begin{equation}\label{3.010}
\prod_{i+j\le 2n}(a_ib_j)^{w(i,j)}.
\notag
\end{equation}
\end{claim}
\begin{proof}
Using claim \ref{Cl3.08} repeatedly we see that
\begin{align}
&\prod_{i+j\le 2n}(a_ib_j)^{w(i,j)}=\prod_{x\in\mathbb{Z},t\le 2n-1}
(a_{(t+x+1)/2}b{(t-x+1)/2})^{\omega(x,t)}\notag\\
&=\prod_{|m|<n}a_{n+m}^{\eta^+(2m,2n-1)}b_{n-m}^{\eta^-(2m,2n-1)}
\prod_{x\in\mathbb{Z},t\le 2n-1}
(a_{(t+x+1)/2}b{(t-x+1)/2})^{T\omega(x,t)}
\notag
\end{align}
Repeated use of this identity toghether with claim \ref{Cl3.05} proves
the claim.
\end{proof}

It is now easy to see that (\ref{0.20}) holds.
\begin{proposition}\label{P3.01}
Set $G(i,j)=h(i-j,i+j-1)$. Then
\begin{equation}\label{3.011}
G(i,j)=\max((g(i-1,j),G(i,j-1))+w(i,j)
\end{equation}
for $i,j\ge 1$.
\end{proposition}
\begin{proof} We have that
\begin{align}
&h(i-j,i+j-1)\notag\\
&=\max(h(i-j-1,i+j-2), h(i-j,i+j-2), h(i-j+1,i+j-2)+w(i,j)\notag\\
&=\max(G(i-1,j),h(i-j,i+j-2), G(i,j-1)+w(i,j).
\notag
\end{align}
Since $h(i-j-1,i+j-2)-h(i-j,i+j-2)=\eta^-(i-j-1,i+j-2)\ge 0$, this last 
expression equals the right hand side of (\ref{3.011}) and we are done.
\end{proof}

If $\eta_r^+, \eta_r^-$ are the jumps for $h_r$ it follows from the 
assignments of weights that $\eta_r^+(2m,2n-1)=u$ has weight $a_{m+n}^u$
and $\eta_r^-(2m,2n-1)=u$ has weight $b_{n-m}^u$, $|m|<n$, $0\le r<n$.
If we think of the weights as labels transported from the $w(i,j)$:s 
we see that if $w(i,j)=0$ for $i>n$ or $j>n$, we have no labels $a_i$ 
with $i>n$ or $b_j$ with $j>n$ and hence $\eta_r^-(2m,2n-1)=0$ if $m<0$
and $\eta_r^+(2m,2n-1)=0$ if $m>0$. Hence all plus-steps thake place to 
the left of the origin and all minus-steps to the right of the origin.
This is the case discussed in \cite{Jo4}. From this consideration and 
(\ref{0.19}) we obtain.

\begin{proposition}\label{P3.1} If $|K|<N$, then
\begin{equation}\label{3.1}
G(N+K,N-K)=h(2K,2N-1).
\end{equation}
Also, if $w(i,j)=0$ for  $|i|>N$ or $j>N$, then for $0\le K<N$,
\begin{equation}\label{3.2}
G(N-K,N)=h(-2K,2N-1)
\end{equation}
and
\begin{equation}\label{3.3}
G(N,N-K)=h(2K,2N-1).
\end{equation}
\end{proposition}

The discussion of the multi-layer extension of the PNG-growth model discussed
above is closely related to the Viennot/matrix-ball construction, \cite{Sa},
\cite{Fu},\cite{Vi}, of the Robinson-Schensted-Knuth (RSK) correspondence.
Let us briefly discuss the relation.
We can think of (\ref{3.03}) and (\ref{3.04}) 
geometrically as follows. From $(x,t)$ to
$(x-1,t+1)$ we draw a line with multiplicity $\eta^+(x,t)$ and from 
$(x,t)$ to
$(x+1,t+1)$ we draw a line with multiplicity $\eta^-(x,t)$. A line
with multiplicity zero means no line. At $(x,t)$ a line with
multiplicity $\eta^+(x+1,t-1)$ and a line with multiplicity
$\eta^-(x-1,t-1)$ meet and we have a collision/annihilation of size
$T\omega(x,t)$ as given by (\ref{3.03}). If $\eta^+(x+1,t-1)\ge 
\eta^-(x-1,t-1)$, then $\eta^+(x+1,t-1)-\eta^-(x-1,t-1)$ plus-lines
survive and we add $\omega(x,t)$ new lines. Similarly in the other
case. This explains (\ref{3.04}). Assume that $w(i,j)=0$ if $i$ or $j$
is $>N$. If $(w(i,j))_{1\le i,j\le N}$ is a permutation matrix this
gives exactly the ``shadow lines'' of the Viennot construction. We
obtain a mapping to a pair of semi-standard Young tableaux $P$ and $Q$
of shape $\lambda$. The number of $m$:s in the first row of $P$ equals
$\eta^-(-(N-m),N+m-1)$, $m=1,\dots,N$ and the number of $m$:s in the
first column of $Q$ equals $\eta^+(-(N-m),N+m-1)$,
$m=1,\dots,N$. Similarly, the same procedure starting with $T\omega$
instead gives the second rows and so on. Using this line of argument
we obtain

\begin{proposition}\label{P3.2}
Let $(w(i,j))_{1\le i,j\le N}$ be given and set $w(i,j)=0$ if $i$ or
  $j$ is $>N$. The RSK-correspondence maps a submatrix $(w(i,j))_{1\le
  i\le M,1\le j\le N}$, $M\le N$ to a pair of semi-standard Young
  tableaux of shape
  $\lambda(M,N)=(\lambda_1(M,N),\lambda_2(M,N),\dots)$. (Similarly, we
  can consider $(w(i,j))_{1\le i\le N,1\le j\le M}$.) Consider the family
  of height curves $h_i$, $0\le i<N$, obtained from the multi-layer
  PNG process using $(w(i,j))$. Then, for $0\le K<N$, $1\le j\le N$,
\begin{equation}\label{3.15}
\lambda_j(N-K,N)=h_{j-1}(-2K,2N-1)+j-1
\end{equation}
and
\begin{equation}\label{3.16}
\lambda_j(N,N-K)=h_{j-1}(2K,2N-1)+j-1.
\end{equation}
\end{proposition}

If we add vertical line segments to the graphs, $x\to h_i(x,2N-1)$,
$0\le i<N$, we obtain $N$ non-intersecting paths with $h_i(-(2N-1),
2N-1)=h_i(2N-1,2N-1)=-i$. Recall that $h_i(x,2N-1)\equiv 1-i$ for $i\ge
N$ so that at most $N$ paths are ``active''. The paths are described by
particle configurations. Let
\begin{equation}
C_{2N-1}(x)=(h_0(x,2N-1),\dots,h_{N-1}(x,2N-1))
\notag
\end{equation}
and
\begin{equation}
C_{2N-1}=(C_{2N-1}(-M+1),\dots,C_{2N-1}(M-1)),
\notag
\end{equation}
where $M=2N-1$. Note that
$C_{2N-1}(-M)=C_{2N-1}(M)=(0,-1,\dots,-N+1)$. Set
\begin{equation}\label{3.17}
\phi_{2j-1,2j}(x,y)=
\begin{cases}(1-a_{j+N})a_{j+N}^{y-x}
   &\text{if $y\ge x$}\\
 0  &\text{if $y<x$,}
\end{cases}
\end{equation}
\begin{equation}\label{3.18}
\phi_{2j,2j+1}(x,y)=
\begin{cases}
0  &\text{if $y>x$,}\\
(1-b_{N-j})b_{N-j}^{x-y}
   &\text{if $y\le x$}
\end{cases}
\end{equation}
for $|j|<N$ with the convention that $0^0=1$. It follows from the
Lindstr\"om-Gessel-Viennot method or from the Karlin-McGregor theorem
that the weight of the non-intersecting path configuration
corresponding to $C_{2N-1}=\bar{x}$, with weights assigned to jumps as
above, equals
\begin{equation}
\left(\prod_{r=-M}^{M-1}
  \det(\phi_{r,r+1}(x_i^r,x_j^{r+1}))_{i,j=1}^N\right) \frac
1{\prod_{j=1}^M(1-a_j)^N(1-b_j)^N}.
\notag
\end{equation}
The way the weights are related to the ``weights'' of the geometric
random variables as described above shows that
\begin{equation}
\mathbb{P}[C_{2N-1}=\bar{x}]=\frac 1{Z_{n,M}}
\prod_{r=-M}^{M-1}
  \det(\phi_{r,r+1}(x_i^r,x_j^{r+1}))_{i,j=1}^n
\notag
\end{equation}
with $Z_{n,M}$ given by (\ref{0.2})
and $n=N$, $M=2N-1$.
Hence we obtain a measure of the form
(\ref{0.1}). We note that 
\begin{equation}\label{3.19}
Z_{n,M}=\frac {\prod_{j=1}^M(1-a_j)^n(1-b_j)^n}{\prod_{i+j\le M}
  (1-a_ib_j)}. 
\end{equation}
We summarize what we have found in the next proposition.

\begin{proposition}\label{P3.3}
Let $h_i$, $i\ge 0$, be the multi-layer PNG process obtained from
geometric random variables with parameters $a_ib_j$ as defined above
and let $\phi_{r,r+1}$ be defined by (\ref{3.17}) and
(\ref{3.18}). Then,
\begin{align}\label{3.20}
&\mathbb{P}[h_{k-1}(r,2N-1)=x_k^r, 1\le i\le n, |r|<M]\\&=
\frac 1{Z_{n,M}}
\prod_{r=-M}^{M-1}
  \det(\phi_{r,r+1}(x_i^r,x_j^{r+1}))_{i,j=1}^n,
\notag
\end{align}
where $Z_{n,M}$ is given by (\ref{3.19}), $x_i^{-M}=x_i^M=1-i$,
$x_1^r>x_2^r>\dots>x_N^r$ for each $r$, $n=N$ and $M=2N-1$.
\end{proposition}

The fact that the probability measure has this form makes it possible
to compute the correlation functions. Set 
\begin{equation}\label{3.21}
f_{2j-1}(z)=(1-a_{j+N})\sum_{m=0}^\infty a_{j+N}^mz^m=
\frac{1-a_{j+N}}{1-a_{j+N}z} 
\end{equation}
and
\begin{equation}\label{3.22}
f_{2j}(z)=(1-b_{N-j})\sum_{m=0}^\infty b_{N-j}^mz^m=
\frac{1-b_{N-j}}{1-b_{N-j}/z} 
\end{equation}
so that (\ref{0.9}) holds. The interpretation of the correlation
functions given by (\ref{0.6}) in this case is that they give the
probability of finding particles at the specified positions. 
We can take $n\ge N$ in (\ref{3.20}), where $n$ is the number of PNG
height curves. All height curves $h_i$ with $i\ge N$ have to be trivial,
i.e. $h_i\equiv -i$ if $i\ge N$. It follows that the probability of a
certain configuration is independent of $n$ for $n\ge N$.
Thus, we can take the kernel $K^{n,M}$, (\ref{0.5}),
with an arbitrary $n$, $n\ge N$ arbitrary and obtain the same value. 
In particular we
can let $n\to\infty$ and use proposition \ref{P0.3}. It is clear that
all the conditions of this theorem are satisfied when $f_r(z)$ is
given by (\ref{3.21}) and (\ref{3.22}). Let $r=2u$, $s=2v$ both be
even, $|u|,|v|<N$. The expression (\ref{0.13}) becomes
\begin{equation}\label{3.23}
G(z,w)=\frac{\prod_{j=u}^{N-1}\left(\frac{1-b_{N-j}}{1-b_{N-j}z}\right) 
\prod_{j=-N+1}^v\left(\frac{1-a_{N+j}}{1-a_{N+j}/w}\right)}
{\prod_{j=-N+1}^u\left(\frac{1-a_{N+j}}{1-a_{N+j}/z}\right)
\prod_{j=v}^{N-1}\left(\frac{1-b_{N-j}}{1-b_{N-j}w}\right) }.
\end{equation}
We summarize our results for the correlation functions in a theorem.

\begin{theorem}\label{T3.4}
Let the multi-layer PNG process be defined using geometric random
variables $w(i,j)$ with parameter $a_ib_j$, $0<a_i,b_j<1$, and let
$G(z,w)$ be given by (\ref{3.23}). Set
\begin{equation}\label{3.24}
\tilde{K}_N(2u,x;2v,y)=\frac 1{(2\pi i)^2}\int_{\gamma_{r_2}}\frac
{dz}z
\int_{\gamma_{r_1}}\frac {dw}w\frac{w^y}{z^x}\frac z{z-w} G(z,w),
\end{equation}
where $\gamma_r$ is the circle with radius $r$ and center at the
origin, $1-\epsilon<r_1<r_2<1+\epsilon$ with $1+\epsilon<\min(1/b_j)$,
$1-\epsilon>\max(a_j)$ and $|u|,|v|<N$,
$x,y\in\mathbb{Z}$. Furthermore, let
\begin{equation}\label{3.25}
\phi_{2u,2v}(x,y)=\frac 1{2\pi}\int_{-\pi}^\pi
e^{i(y-x)\theta}G(e^{i\theta},e^{i\theta}) d\theta,
\end{equation}
for $u<v$ and $\phi_{2u,2v}(x,y)=0$ for $u\ge v$. Set
\begin{equation}\label{3.26}
K_N(2u,x;2v,y)=\tilde{K}_N(2u,x;2v,y)-\phi_{2u,2v}(x,y).
\end{equation}
Then,
\begin{align}\label{3.27}
&\mathbb{P}[(2u,x_j^{2u})\in\{(2t,h_i(2t,2N-1))\,;\, \text{$|t|<N$ , 
$0\le i<N$}\},
|u|<N, 1\le j\le k_u] \\
&=\det(K_N(2u,x_i^{2u};2v,x_j^{2v}))_{|u|,|v|<N,1\le i\le k_u, 1\le
  j\le k_v}
\notag
\end{align}
for any $x_j^{2u}\in\mathbb{Z}$ and any $k_u\in\{0,\dots,N\}$.
\end{theorem}

Consider the finite-dimensional distribution of $h_0(x,t)=h(x,t)$, the top 
curve,
\begin{equation}\label{3.28}
\mathbb{P}_{n,M}[h_0(2s_i,2N-1)\le\ell_i, 1\le i\le m],
\end{equation}
where $n$ is the number of paths, $M=2N-1$, $n\ge N$, $|s_i|<N$ 
and $\ell_i>-N$. This can also be written
\begin{equation}
\mathbb{P}_{n,M}[\text{no particles in $\{2s_i\}\times (\ell_i,\infty)$}, 
1\le i\le m].
\notag
\end{equation}
This probability is independent of $n\ge N$ and hence we can let $n\to\infty$.
Let $g(2s_i,x)=-\chi_{\ell_i,\infty)}(x)$, $1\le i\le m$, and 
$g(r,x)\equiv 0$ if $r$ is not equal to one of the $2s_i$. Hence, by 
proposition \ref{P2.01},
\begin{equation}\label{3.32}
\mathbb{P}_{N,M}[h_0(2s_i,2N-1)\le\ell_i, 1\le i\le m]
=\det (I+gK_N)_{L^2(\Lambda_M)}.
\end{equation}
This formula can be used to study the convergence in distribution
of the rescaled height curve.

\section{Asymptotics}
We will consider the asymptotics of the kernel (\ref{3.26}) in the
case $a_i=b_i=\alpha$ for all $i\ge 1$, so that $w(i,j)$ are geometric
random variables with parameter $q=\alpha^2$. The function $G(z,w)$ in
(\ref{3.23}) then becomes
\begin{equation}\label{4.1}
G(z,w)=(1-\alpha)^{2(v-u)}\frac{(1-\alpha/z)^{N+u}}{(1-\alpha z)^{N-u}}
\frac{(1-\alpha w )^{N-v}}{(1-\alpha/w)^{N+v}}.
\end{equation}
Write
\begin{equation}
F_{u,x}(z)=\frac 1{z^{x+N+u}}\frac{(z-\alpha)^{N+u}}{(1-\alpha
  z)^{N-u}}, 
\notag
\end{equation}
so that, by (\ref{3.24}),
\begin{equation}\label{4.2}
\tilde{K}_N(2u,x;2v,y)=\frac {(1-\alpha)^{2(v-u)}}{(2\pi
  i)^2}\int_{\gamma_{r_2}} \frac{dz}z\int_{\gamma_{r_1}} \frac{dw}w
  \frac{z}{z-w} F_{u,x}(z)F_{-v,y}(\frac 1w),
\end{equation}
where $\alpha<r_1<r_2<1/\alpha$.

Set $\mu=m/N$, $\mu'=m'/N$, $\beta=u/N$, $\beta'=-v/N$ and
\begin{align}
f_{\mu,\beta}(z)&=\frac 1N \log F_{u,m-N}(z)\notag\\&=(1+\beta)\log
(z-\alpha)-(1-\beta)\log (1-\alpha z)-(\mu+\beta)\log z.
\notag
\end{align}
Then, $f'(z)=P(z)/Q(z)$, where $Q(z)=z(z-\alpha)(1-\alpha z)$ and
\begin{equation}
P(z)=\alpha (\mu-\beta)[z^2+\frac
{1-\alpha^2-\mu(1+\alpha^2)}{\alpha(\mu-\beta)}
z+\frac{\mu+\beta}{\mu-\beta}].
\notag
\end{equation}
We will write
\begin{equation}\label{4.1'}
p=p(\mu,\beta)=\frac{\mu(1+\alpha^2)-(1-\alpha^2)}{2\alpha(\mu-\beta)}\quad
;\quad q=a(\mu,\beta)=\frac{\mu+\beta}{\mu-\beta}.
\end{equation}
The critical points of $f$ are $z_c^{\pm}=p\pm\sqrt{p^2-q}$ and we
obtain a double critical point if $p^2=q$ which gives
\begin{equation}\label{4.2'}
\mu=\mu_c(\beta)=\frac{1+\alpha^2}{1-\alpha^2}
+\sqrt{\left(\frac{1+\alpha^2}{1-\alpha^2}\right)^2
  -1-\frac{4\alpha^2\beta^2}{(1-\alpha^2)^2} }
\end{equation}
and
\begin{equation}\label{4.3}
p_c=p(\mu_c,\beta)=\frac{2\alpha
  +(1+\alpha^2)\sqrt{1-\beta^2}}{1+\alpha^2+
2\alpha\sqrt{1-\beta^2}-\beta(1-\alpha^2)}. 
\end{equation}
Set
\begin{equation}\label{4.4}
d=\frac{\alpha^{1/3}(1+\alpha)^{1/3}}{1-\alpha} \quad,\quad
d'=\frac{1-\alpha}{1+\alpha} d.
\end{equation}
It will be convenient to write
\begin{equation}\label{4.5}
u=\frac 1{d'}\tau N^{2/3}\quad,\quad
v=\frac 1{d'}\tau' N^{2/3},
\end{equation}
since $N^{2/3}$ is the right scale for $u$ and $v$ if we want a
non-trivial limit. The correct way of writing $x$ and $y$ will turn
out to be
\begin{equation}\label{4.6}
x=N(\mu_c(\beta)-1)+\xi dN^{1/3}\quad,\quad
y=N(\mu_c(\beta')-1)+\xi' dN^{1/3}.
\end{equation}
We will assume that $|\tau|, |\tau'|, |\xi|, |\xi'|$ are $\le \log N$.

The paths of integration can be deformed into
\begin{align}\label{4.7}
&\Gamma:\mathbb{R}\ni t'\to z(t')=p_c(\beta)+\frac{\eta}{dN^{1/3}}-
\frac{it'}{dN^{1/3}}\doteq p-it,\\
&\Gamma':\mathbb{R}\ni s'\to w(s')=p_c(\beta')+\frac{\eta'}{dN^{1/3}}-
\frac{is'}{dN^{1/3}}\doteq (p'-it)^{-1},
\end{align}
where $\eta,\eta'>0$ will be appropriately chosen; we will require that
\begin{equation}\label{4.8}
(p_c(\beta)+\frac{\eta}{dN^{1/3}})(p_c(\beta')+\frac{\eta'}{dN^{1/3}})>1.
\end{equation}
In that case we have, by Cauchy's theorem,
\begin{equation}\label{4.9}
\tilde{K}_N(2u,x;2v,y)=-\frac {(1-\alpha)^{2(v-u)}}{(2\pi
  i)^2}\int_{\Gamma} \frac{dz}z\int_{\Gamma'} \frac{dw}w
  \frac{z}{z-w} F_{u,x}(z)F_{-v,y}(\frac 1w).
\end{equation}
We first estimate this integral and then we will compute its
asymptotics using a saddle-point argument. Choose $\mu$ so that
$p(\mu,\beta)=p\doteq p_c(\beta)+\eta/dN^{1/3}$ as in (\ref{4.7}), and
let $q=q(\mu,\beta)$ be the value we get with this $\mu$. This is
possible by formula (\ref{4.1'}) with $\mu>\mu_c$. We can write
\begin{equation}
F_{u,x}(z)=\frac 1{z^{x+N-\mu N}}\frac 1{z^{(\mu+\beta)N}}\frac
{(z-\alpha)^{N+u}}{(1-\alpha z)^{N-u}},
\notag
\end{equation}
and then let $z=p-it$, $t\in\mathbb{R}$, and take the absolute value
to get
\begin{equation}\label{4.10}
|F_{u,x}(p-it)|^2=\frac 1{(p^2+t^2)^{x+N-\mu N}}e^{2Nh(t)},
\end{equation}
where
\begin{equation}
2h(t)=(1+\beta)\log A-(1-\beta)\log B-(\mu+\beta)\log C,
\notag
\end{equation}
with
\begin{equation}
A=(p-\alpha)^2+t^2=1-2\alpha p+\alpha^2+\frac{2\beta}{\mu-\beta}+p^2-q+t^2
\notag
\end{equation}
\begin{equation}
B=(1-\alpha p)^2+\alpha^2t^2=1-2\alpha
p+\alpha^2+\frac{2\alpha^2\beta}
{\mu-\beta}+\alpha^2(p^2-q+t^2)
\notag
\end{equation}
\begin{equation}
C=p^2+t^2=1+\frac{2\beta}{\mu-\beta}+p^2-q+t^2.
\notag
\end{equation}
Note that 
\begin{equation}
(\mu-\beta)A=(1-\alpha^2)\beta+1-\alpha^2+(\mu-\beta)(p^2-q+t^2)
\notag
\end{equation}
\begin{equation}
(\mu-\beta)B=1-\alpha^2-(1-\alpha^2)\beta+\alpha^2(\mu-\beta)(p^2-q+t^2)
\notag
\end{equation}
\begin{equation}
(\mu-\beta)C=(\mu-\beta)(p^2-q+t^2)+\mu+\beta.
\notag
\end{equation}
A computation now gives
\begin{align}\label{4.11}
h'(t)=&\frac{t(p^2-q+t^2)}{(\mu-\beta)ABC}
\{(1-\alpha^2)^2-(1-\alpha^4)(\mu+\beta)+
(1+\alpha^4)\mu\beta+2\alpha^2\beta^2\\&-\alpha^2(\mu-\beta)^2(p^2-q+t^2)\}. 
\notag
\end{align}
Another computation shows that $p^2-q\ge 0$. To leading order we have
$p^2-q\approx 2\eta/dN^{1/3}$. Recall that $t=t'/dN^{1/3}$. When $t$
is large we have $h'(t)\approx -(\mu-\beta)/t$, and we also have the
estimate
\begin{equation}\label{4.12}
h'(t)\le\frac{t}{(\mu-\beta)ABC}(p^2-q)
\{(1-\alpha^2)^2-(1-\alpha^4)(\mu+\beta)+
(1+\alpha^4)\mu\beta+2\alpha^2\beta^2\}. 
\end{equation}
for $t\ge 0$. Consider the case $0\le t'\le N^\gamma$; the case
$-N^\gamma\le t'\le 0$ is analogous by symmetry. Here
$0<\gamma<1/3$. Using (\ref{4.12}) we see that $h(t)-h(0)\approx
-2\eta {t'}^2/N$ and we can show that
\begin{equation}\label{4.13}
Nh'(t)\le -\frac 32\eta {t'}^2+Nh(0)
\end{equation}
for $|t'|\le N^\gamma$, $0<\gamma<1/3$, and $N$ sufficiently
large. Define $h_\ast(t,p)$ by
\begin{equation}
e^{2Nh_\ast(t,p)}=\frac 1{(p^2+t^2)^{(\mu_c+\beta)N}}
\frac{[(p-\alpha)^2+t^2]^{(1+\beta)N}}{[(1-p\alpha)^2+
\alpha^2t^2]^{(1-\beta)N}}.
\notag
\end{equation}
A computation, compare (\ref{4.20}) below, gives
\begin{equation}
e^{h_\ast(0,p_c)}\sim (1-\alpha)^{2u}e^{\frac 13 \tau^3}
\notag
\end{equation}
and $h_\ast(0,p)=h_\ast(0,p_c)+d^3(p-p_c)^3/3+ \dots=
h_\ast(0,p_c)+\eta^3/3N+\dots$. Consequently,
\begin{equation}
e^{Nh(0)}=\frac 1{p^{N(\mu-\mu_c)}}e^{Nh_\ast(t,p)}\sim
\frac {(1-\alpha)^{2u}}{p^{N(\mu-\mu_c)}} e^{(\tau^3+\eta^3)/3}.
\notag
\end{equation}
Combining this with (\ref{4.13}) gives
\begin{equation}
|F_{u,x}(p-it)|\le \frac C{(p^2+t^2)^{(x+N(1-\mu))/2}}
\frac {(1-\alpha)^{2u}}{p^{N(\mu-\mu_c)}} e^{(\tau^3+\eta^3)/3-
3\eta {t'}^2/2}.
\notag
\end{equation}
Write
\begin{equation}
\frac 1{(p^2+t^2)^{(x+N(1-\mu))/2}}=
\left|\frac 1{(p-it)^{x+N(1-\mu_c)}}\right|\left(\frac
  {p^2+t^2}{p^2}\right)^{N(\mu-\mu_c)/2}. 
\notag
\end{equation}
Further computation shows that 
\begin{equation}
\left|\frac 1{(p-it)^{x+N(1-\mu_c)}}\right|\sim e^{-\xi\tau-\xi\eta}
\notag
\end{equation}
\begin{equation}
(1+t^2/p^2)^{N(\mu-\mu_c)/2}\sim e^{\eta {t'}^2}.
\notag
\end{equation}
Collecting the estimates we find
\begin{equation}\label{4.14}
|F_{u,x}(p_c(\beta)+\frac \eta{dN^{1/3}}-\frac{it'}{dN^{1/3}})|\le
C(1-\alpha)^{2u} e^{\frac 13(\tau^3+\eta^3)-\xi\tau-\xi\eta-\frac
  {\eta}2{t'}^2} 
\end{equation}
for $|t'|\le N^{\gamma}$, $0<\eta<N^\gamma$, $0<\gamma<1/3$.

Using (\ref{4.11}), (\ref{4.12}) and the other estimates above we see
that the contribution to the integral from $|t'|\ge N^\gamma$ and/or
$|s'|\ge N^{\gamma}$ is $\le C\exp(-cN^{2\gamma})$ for some constant
$c>0$. Hence, using the parametrization (\ref{4.7}) in (\ref{4.9}) we
can restrict to $|t'|\le N^\gamma$, $|s'|\le N^\gamma$. We can 
use (\ref{4.14}) if we want an estimate of the integral. To get the
asymptotics we make a local saddle-point argument. 

To leading order we have $p_c(\beta)=1+\tau/dN^{1/3}$,
$p_c(\beta')=1-\tau'/dN^{1/3}$ and hence the condition (\ref{4.8})
requires
\begin{equation}\label{4.15}
\tau-\tau'+\eta+\eta'>0.
\end{equation}
We will use the parametrizations (\ref{4.7}) and consider the integral
\begin{align}\label{4.16}
-\frac {(1-\alpha)^{2(v-u)}} {(2\pi i)^2} \int_{|t|\le N^\gamma}dt
\int_{|s|\le N^\gamma}ds&
\frac{z'(t)}{z(t)}\frac{w'(t)}{w(t)}\frac{z(t)}{z(t)-w(s)}
\frac{w(s)^{y+N(1-\mu_c(\beta'))}}{z(t)^{x+N(1-\mu_c(\beta))}}\\&\times
e^{Nf_{\mu_c(\beta),\beta}(z(t)) +Nf_{\mu_c(\beta'),\beta'}(1/w(s))}.
\notag\end{align}
Now,
\begin{align}\label{4.17}
&Nf_{\mu_c(\beta),\beta}(p_c(\beta)+\frac 1{dN^{1/3}}(\eta-it))\\&=
Nf_{\mu_c(\beta),\beta}(p_c(\beta))+\frac i3\left(\frac 1{2d^3}
f_{\mu_c(\beta),\beta}^{(3)}(p_c(\beta))\right)(t+i\eta)^3+r_N(t)\notag\\
&=Nf_{\mu_c(\beta),\beta}(p_c(\beta))+\frac i3(t+i\eta)^3+r_N(t)
\notag
\end{align}
where the remainder term $r_N(t)$ can be neglected for $|t|\le
N^\gamma$. Also, $z'(t)=-i/dN^{1/3}$, $w'(s)=iw(s)^2/dN^{1/3}$ and
\begin{equation}\label{4.18}
\frac{w(s)}{z(t)-w(s)}\sim -\frac{dN^{1/3}}{\tau'-\tau+i(t+i\eta+s+i\eta')}.
\end{equation}
Furthermore
\begin{equation}\label{4.19}
\frac{w(s)^{y+N(1-\mu_c(\beta'))}}{z(t)^{x+N(1-\mu_c(\beta))}}\sim 
e^{\xi'\tau'-\xi\tau+i\xi(t+i\eta)i\xi'(s+i\eta')}.
\end{equation}
We also need to compute
\begin{equation}
e^{Nf_{\mu_c(\beta),\beta}(p_c(\beta))}=\frac
{(p_c(\beta)-\alpha)^{N+u}}{(1-\alpha p_c(\beta))^{N-u}}\frac
1{p_c^{N\mu_c(\beta)-u}} .
\notag
\end{equation}
Using the formulas (\ref{4.2'}) and (\ref{4.3}) above a rather long
computation, which we omit, shows that
\begin{equation}\label{4.20}
e^{Nf_{\mu_c(\beta),\beta}(p_c(\beta))}\sim (1-\alpha)^{2u}e^{\tau^3/3}.
\end{equation}
Inserting (\ref{4.17}) - (\ref{4.20}) into (\ref{4.16}) we see that,
provided (\ref{4.15}) holds,
\begin{align}\label{4.21}
\lim_{N\to\infty} dN^{1/3}\tilde{K}_N&(2\frac {1+\alpha}{1-\alpha}\frac
1d N^{2/3}\tau,\frac {2\alpha}{1-\alpha}N +(\xi-\tau^2)dN^{1/3};\\&
2\frac {1+\alpha}{1-\alpha}\frac
1d N^{2/3}\tau',\frac {2\alpha}{1-\alpha}N +(\xi'-{\tau'}^2)dN^{1/3})
\notag\\
&=-\frac 1{4\pi^2}e^{\frac 13(\tau^3-{\tau'}^3)+\xi'\tau'-\xi\tau}
\int_{\im z=\eta}\int_{\im w=\eta'}\frac {e^{i\xi z+i\xi'w+\frac i3
    (z^3+w^3)}}{\tau'-\tau+i(z+w)} dzdw.
\notag
\end{align}
Here we have used $\mu(\beta_c)=\frac
{1+\alpha}{1-\alpha}-\frac{\alpha}{1-\alpha^2}\beta^2+\dots$. We also
want to compute the corresponding limit of (\ref{3.25}) with $G(w,w)$ given
by (\ref{4.1}), i.e. we consider, $u<v$, 
\begin{equation}\label{4.22}
\phi_{2u,2v}(x,y)=\frac{(1-\alpha)^{2(v-u)}}{2\pi }\int_{-\pi}^\pi
e^{i(y-x)\theta+(v-u)\log (1+\alpha^2-2\alpha\cos\theta)}d\theta.
\end{equation}
If we set $g(\theta)=\log (1+\alpha^2-2\alpha\cos\theta)$, then
\begin{equation}
g'(\theta)=\frac{2\alpha\sin\theta}{1+\alpha^2-2\alpha\cos\theta},
\notag
\end{equation}
and we see that $g(\theta)$ has a quadratic minimum at
$\theta=0$. Hence, we can immediately both compute the asympotics of
and estimate the integral in (\ref{4.22}) when
$x=2\alpha(1-\alpha)^{-1} N+(\xi-\tau^2)dN^{1/3}$, 
$y=2\alpha(1-\alpha)^{-1} N+(\xi'-{\tau'}^2)dN^{1/3}$,
$u=\frac{1+\alpha}{1-\alpha}d^{-1}N^{2/3}\tau$ and 
$v=\frac{1+\alpha}{1-\alpha}d^{-1}N^{2/3}\tau'$. We obtain
\begin{align}\label{4.23}
\lim_{N\to\infty} dN^{1/3}\phi_{2u,2v}(x,y)&=
\frac 1{2\pi}\int_{\mathbb{R}}
e^{i(\xi'-\xi+\tau^2-{\tau'}^2)t-(\tau'-\tau)t^2} dt\\
&=\frac 1{\sqrt{4\pi
    (\tau'-\tau)}}e^{-(\xi'-\xi+\tau^2-{\tau'}^2)^2/(\tau'-\tau)}.
\notag 
\end{align}
We want to identify the right hand side of (\ref{4.21}) combined with
(\ref{4.23}) with the extended Airy kernel. This can be done using the
proposition \ref{L4.1}.
Combining this double integral formula for the
extended Airy kernel with (\ref{4.21}) and (\ref{4.23}) we obtain
the following result.

\begin{proposition}\label{P4.2}
Let $d=(1-\alpha)^{-1}\alpha^{1/3}(1+\alpha)^{1/3}$ and let $K_N$ be
given by (\ref{3.26}). Then
\begin{align}\label{4.27}
\lim_{N\to\infty} dN^{1/3}K_N&(2\frac
{1+\alpha}{1-\alpha}d^{-1}N^{2/3}\tau, 
\frac{2\alpha}{1-\alpha}N+(\xi-\tau^2)dN^{2/3};\notag\\
&2\frac
{1+\alpha}{1-\alpha}d^{-1}N^{2/3}\tau', 
\frac{2\alpha}{1-\alpha}N+(\xi'-{\tau'}^2)dN^{2/3})\notag\\
&=e^{(\tau^3-{\tau'}^3)/3+\xi'\tau'-\xi\tau}A(\tau,\xi;\tau'\xi')
\end{align}
uniformly for $\xi,\xi',\tau,\tau'$ in a compact set.
\end{proposition}

We can now combine the formula \ref{3.32}, theorem \ref{T3.4},
proposition \ref{P4.2} and some estimates of $K_N$, which can be
obtained from the asymptotic analysis above, to prove the
following theorem on convergence in distribution to the Airy
process. A complete proof requires a control of the convergence of the
Fredholm expansions but we will not present the details. The
individual determinants in the Fredholm expansion can be estimated
using the Hadamard inequality. Compare with theorem 1.2 and lemma 3.1 
in \cite{Jo1}.

\section{A functional limit theorem}
\subsection{A moment estimate}
Consider the PNG height functions $h_k(x,2N-1)$ defined in sect. 3. Set
\begin{equation}
t_j=\frac{j}{cN^{2/3}}
\notag
\end{equation}
where $c=(1+\alpha)(1-\alpha)^{-1}d^{-1}$, $j\in\mathbb{Z}$. The
normalized height functions are
\begin{equation}
H_{N,k}(t_j)=\frac 1{dN^{1/3}}(h_k(2j,2N-1)-\frac
{2\alpha}{1-\alpha} N),
\notag
\end{equation}
with $d$ as in (\ref{4.4}), $k\in\mathbb{N}$. For a given function
$f:\mathbb{R}\to\mathbb{C}$ we write
\begin{equation}
f_N(x)=f(\frac 1{dN^{1/3}} (x-\frac
{2\alpha}{1-\alpha} N)).
\notag
\end{equation}
Assume that there is a $K$ such that $f(x)=0$ for $x\le K$. Define
\begin{equation}\label{5.1}
H_N(f,t_j)=\sum_{k=0}^\infty f_N(h_k(2j,2N-1))=\sum_{k=0}^\infty
f(H_{N,k}(t_j)).
\end{equation}

\begin{lemma}\label{L5.1} 
Assume that $f$ is a $C^\infty$ function and that
  there are constants $K_1$ and $K_2$ such that $f(x)=0$ if $x\le K_1$,
  and $f(x)$ equals a constant if $x\ge K_2$. There is a constant
  $C(f,\alpha)$ so that
\begin{equation}\label{5.2}
\mathbb{E}[(H_N(f,t_u)-H_N(f,t_v))^4]\le
C(f,\alpha)e^{-|t_u|^3}|t_u-t_v|^2, 
\end{equation}
for $|t_u-t_v|\le 1$ and $|t_u|,|t_v|\le\log N$.
\end{lemma}

\begin{proof}
The proof is rather long and complicated. We will outline the main
ideas and steps in the argument without giving full details. The left
hand side of (\ref{5.2}) can be written
\begin{equation}\label{5.3}
\sum_{k_1,k_2,k_3,k_4=1}^\infty
\mathbb{E}[\prod_{r=1}^4(f_N(h_{k_r}(2u,2N-1))-f_N(h_{k_r}(2v,2N-1)))].
\end{equation}
We can rewrite this using formula (\ref{3.27}) in theorem \ref{T3.4}. Let
us write the kernel $K_N(2u,x;2v,y)$ in (\ref{3.26}) as
$K_{uv}(x,y)$. 
We will use the following notation:
\begin{equation}\label{1.0}
K^{r_1\,\, r_2\dots r_m}_{s_1\,\, s_2\dots s_m}
\left(\begin{matrix} x_1 & x_2 &\dots &x_m \\ y_1 & y_2 & \dots &
  y_m\end{matrix}\right) =\det (K(r_i,x_i;s_j,y_j))_{i,j=1}^m,
\end{equation}
and we will also write
\begin{equation}\label{1.1}
K(\begin{matrix} r_1,x_1 & r_2,x_2 &\dots & r_m,x_m
\end{matrix}) =\det (K(r_i,x_i;r_j,x_j))_{i,j=1}^m.
\end{equation}
Furthermore, we will write
\begin{equation}\label{5.3'}
D_{u_1,\dots,u_m}(x_1,\dots,x_m)=K(2u_1, x_1\,\,\, 2u_2, x_2 \dots
2u_m,x_m).
\end{equation}
Set
$h_N(x_1,x_2,x_3)=-6[f_N(x_1)^2f_N(x_2)f_N(x_3)+
f_N(x_1)f_N(x_2)^2f_N(x_3)$\linebreak 
$-f_N(x_1)f_N(x_2)f_N(x_3)]$, which is
synmmetric under permutation of $x_1$ and $x_2$. Then, the sum in
(\ref{5.3}) can be written
\begin{align}\label{5.4}
&\sum_{x\in\mathbb{Z}^4}
f_N(x_1)f_N(x_2)f_N(x_3)f_N(x_4)[D_{uuuu}(x_1,x_2,x_3,x_4) -
4D_{uuuv}(x_1,x_2,x_3,x_4)\\&+6D_{uuvv}(x_1,x_2,x_3,x_4)
-4D_{uvvv}(x_1,x_2,x_3,x_4)+D_{vvvv}(x_1,x_2,x_3,x_4)]\notag\\
&+\sum_{x\in\mathbb{Z}^3}\{6f_N(x_1)^2f_N(x_2)f_N(x_3)[D_{uuu}(x_1,x_2,x_3)+
D_{vvv}(x_1,x_2,x_3)]\notag\\
&+h_N(x_1,x_2,x_3)D_{uuv}(x_1,x_2,x_3)+
h_N(x_3,x_2,x_1)D_{uvv}(x_1,x_2,x_3)\} 
\notag\\
&+\sum_{x\in\mathbb{Z}^2} 2(f_N(x_1)^3f_N(x_2)+f_N(x_1)f_N(x_3)^3)
[D_{uu}(x_1,x_2)-2D_{uv}(x_1,x_2)+D_{vv}(x_1,x_2)]
\notag\\
&+\sum_{x\in\mathbb{Z}^2}3f_N(x_1)^2f_N(x_2)^2[D_{uu}(x_1,x_2)+
2D_{uv}(x_1,x_2)+D_{vv}(x_1,x_2)]
\notag\\
&+\sum_{x_1\in\mathbb{Z}} 2f_N(x_1)^4[K_{uu}(x_1,x_1)+K_{vv}(x_1,x_1)]
\notag\\
&\doteq \Sigma_1+\Sigma_2+\Sigma_3+\Sigma_4+\Sigma_5.
\notag
\end{align}
We have $K_{uv}=\tilde{K}_{uv}-\phi$ if $u<v$ and
$K_{uv}=\tilde{K}_{uv}$ if $u\ge v$. Here we have written
$\phi=\phi_{u,v}$. Set
$\Delta K_{uv}=\tilde{K}_{uv}-\tilde{K}_{uu}$,
$\Delta K_{vu}=\tilde{K}_{vu}-\tilde{K}_{uu}$
$\Delta K_{vv}=\tilde{K}_{uv}-\tilde{K}_{uu}$. We see from (\ref{4.23}) 
that $\phi$ acts like a kind of approximate $\delta$-function. This
will be important for the cancellation between different terms in
(\ref{5.4}). The argument goes as follows. We will take out all terms
in (\ref{5.4}) containing $\phi$ and combine them with other terms so
that we get cancellation. We will then expand in $\Delta K_{uv}$,
$\Delta K_{vu}$ and $\Delta K_{vv}$. The terms linear in $\Delta K$
will cancel and what will remain will be terms containg $(\Delta K)^2$
or higher powers. They will give a contribution proportional to
$|t_u-t_v|^2$ which is what we want. 

In the computations below we use symmetries and also relabelling of
variables. Expand in $\phi$ and in the terms linear in $\phi$ we
expand in $\Delta K$. 
Let $\tilde{D}$ denote the same
object as in (\ref{5.3'}) but with $K$ replaced by $\tilde{K}$.
We find
\begin{align}\label{5.5} 
\Sigma_1&=\sum_x f_N(x_1)f_N(x_2)f_N(x_3)f_N(x_4)[\tilde{D}_
{uuuu}(x_1,x_2,x_3,x_4) -
4\tilde{D}_{uuuv}(x_1,x_2,x_3,x_4)\\&+6\tilde{D}_{uuvv}(x_1,x_2,x_3,x_4)
-4\tilde{D}_{uvvv}(x_1,x_2,x_3,x_4)+
\tilde{D}_{vvvv}(x_1,x_2,x_3,x_4)]\notag\\
&+\sum_x 24\phi(x_1,x_4)K_{uu}^{uu}
\left(\begin{matrix} x_3 & x_4 \\ x_1 & x_2
  \end{matrix}\right) [\Delta K_{vu}(x_2,x_3)+\Delta
K_{uv}(x_2,x_3)-\notag\\
&\Delta K_{vv}(x_2,x_3)]f_N(x_1)f_N(x_2)f_N(x_3)f_N(x_4)
\notag\\
&-\sum_x 12\phi(x_1,x_4)\phi_(x_2,x_3)\tilde{K}_{uu}^{vv}
\left(\begin{matrix} x_3 & x_4 \\ x_1 & x_2
  \end{matrix}\right)f_N(x_1)f_N(x_2)f_N(x_3)f_N(x_4)
\notag\\
&\sum_x (\text{terms with $\Delta K^2$}).
\notag
\end{align}
We will give a brief discussion of the $\Delta K^2$ terms below. Also
we will see then that terms containing
\begin{equation}\label{5.6}
\Delta K_{vu}(x,y)+\Delta K_{uv}(x,y)-\Delta K_{vv}(x,y)
\end{equation}
will give a contribution proportional to $|t_u-t_v|^2$. If we expand
the $\tilde{D}$-part of (\ref{5.5}) in $\Delta K$ we will see that the
terms linear in $\Delta K$ cancel out. Since obviously the 0:th order
term equals zero we are left with $\Delta K^2$-terms. The term
containg two $\phi$-factors will be combined with other terms below.

We expand $\Sigma_2$ similarly. The part linear in $\Delta K$ is
\begin{align}\label{5.7}
-\sum_x 24f_N(x_1)^2f_N(x_2)f_N(x_3)&\tilde{K}_{uu}^{uu}
\left(\begin{matrix} x_3 & x_1 \\ x_1 & x_2
  \end{matrix}\right)
[\Delta K_{vu}(x_2,x_3)\\&+\Delta
K_{uv}(x_2,x_3)-\Delta
K_{vv}(x_2,x_3)].
\notag
\end{align}
Actually this sum can be combined with the corresponding term in
(\ref{5.5}) to get some cancellation, see the $\phi$-calculations
below, but we can also use the fact that (\ref{5.6}) has the right
order. We also get $\Delta K^2$-terms and a term linear in $\phi$,
\begin{align}\label{5.8} 
&\sum_x 12\phi(x_1,x_3)[f_N(x_1)^2f_N(x_2)f_N(x_3)
-f_N(x_1)f_N(x_2)f_N(x_3)^2]\\
&\times \left\{\tilde{K}_{uu}^{uv}
\left(\begin{matrix} x_2 & x_3 \\ x_1 & x_2
  \end{matrix}\right)-
\tilde{K}_{uv}^{vv}
\left(\begin{matrix} x_2 & x_3 \\ x_1 & x_2
  \end{matrix}\right)\right\}
\notag\\
&+\sum_x 12\phi(x_1,x_3)f_N(x_1)f_N(x_2)^2f_N(x_3)\left\{\tilde{K}_{uu}^{uv}
\left(\begin{matrix} x_2 & x_3 \\ x_1 & x_2
  \end{matrix}\right)+
\tilde{K}_{uv}^{vv}
\left(\begin{matrix} x_2 & x_3 \\ x_1 & x_2
  \end{matrix}\right)\right\}
\notag\\
&\doteq \Sigma_a+\Sigma_b.
\notag
\end{align}

Consider next $\Sigma_3$. We get a term linear in $\phi$,
\begin{equation}\label{5.9}
-\sum_x
 4\phi(x_1,x_2)K_{vu}(x_2,x_1)[f_N(x_1)f_N(x_2)^3+f_N(x_1)^3f_N(x_2)] 
\end{equation}
a term linear in $\Delta K$,
\begin{align}\label{5.10}
\sum_x
 4(f_N(x_1)f_N(x_2)^3+f_N(x_1)^3&f_N(x_2)])[\Delta K_{vu}(x_1,x_2)+\\&\Delta
K_{uv}(x_1,x_2)-\Delta
K_{vv}(x_1,x_2)]K_{uu}(x_2,x_1),
\notag
\end{align}
and $\Delta K^2$-terms. In (\ref{5.10}) we again have the expression
(\ref{5.6}).

The leading term in $\Sigma_4$ is
\begin{equation}\label{5.11}
\sum_x 3f_N(x_1)^2f_N(x_2)^2\left[\tilde{K}_{uu}^{uu}
\left(\begin{matrix} x_1 & x_2 \\ x_1 & x_2
  \end{matrix}\right)+2
\tilde{K}_{uv}^{uv}
\left(\begin{matrix} x_1 & x_2 \\ x_1 & x_2
  \end{matrix}\right)+
\tilde{K}_{vv}^{vv}
\left(\begin{matrix} x_1 & x_2 \\ x_1 & x_2
  \end{matrix}\right)
\right]
\end{equation}
and we also have a term linear in $\phi$,
\begin{equation}\label{5.12}
\sum_x 6f_N(x_1)^2f_N(x_2)^2\phi(x_1,x_2) \tilde{K}_{vu}(x_2,x_1).
\end{equation}
Finally we have $\Sigma_5$ which is 
\begin{equation}\label{5.13}
\sum_{x_1} 2f_N(x_1)^4[\tilde{K}_{uu}(x_1,x_1)+\tilde{K}(x_1,x_1)].
\end{equation}
 When calculating the cancellations involving the $\phi$-terms we will
 combine the double $\phi$-term in (\ref{5.5}) with $\Sigma_b$ in
 (\ref{5.8}) and (\ref{5.11}). Also we will combine (\ref{5.9}),
 (\ref{5.12}) and (\ref{5.13}). We will discuss this second case first
 in some detail and then the first case more briefly. The term
 $\Sigma_a$  is similar and finally we will indicate what is involved
 in estimating (\ref{5.6}) and the $\Delta K^2$-terms.

We want to estimate
\begin{align}\label{5.14}
&\sum_{x,y\in\mathbb{Z}} \phi(x,y)
\bar{K}_{vu}(y,x)[6f_N(x_1)^2f_N(x_2)^2- 4f_N(x_1)f_N(x_2)^3
\\
&-4f_N(x_1)^3f_N(x_2)]+\sum_{y\in\mathbb{Z}} (\tilde{K}_{uu}(y,y)+
\tilde{K}_{vv}(y,y))f_N(y)^4.
\notag
\end{align}
Here we have made a symmetrization in $x$ and $y$ by setting
\begin{equation}
\bar{K}_{vu}(x,y)=\frac 12[\tilde{K}_{vu}(x,y)+\tilde{K}_{vu}(y,x)];
\notag
\end{equation}
note that $\phi(x,y)$ is symmetric in $x$ and $y$.
Next, we will introduce some notation and some formulas that will be
used. Set
\begin{equation}\label{5.141}
g(z)=-\frac{\alpha}{(1-\alpha)^2} (z+\frac 1z -2)
\end{equation}
and
\begin{equation}\label{5.15}
G_{ab}^\ast (z,w)=\left(\frac{1-\alpha/z}{1-\alpha z}\right)^N 
\left(\frac{1-\alpha w}{1-\alpha/w}\right)^N (1+g(z))^a(1+g(w))^b
\frac 1{w(z-w)}
\end{equation}
so that
\begin{equation}\label{5.150}
\tilde{K}_{ab}(x,y)=\frac 1{(2\pi i)^2}\int_{\gamma_{r_2}} dz
\int_{\gamma_{r_2}}dw  G_{ab}^\ast (z,w) \frac {w^y}{z^x},
\end{equation}
and
\begin{equation}\label{5.151}
\phi(x,y)=\frac 1{2\pi} \int_{-\pi}^\pi e^{i(y-x)\theta}
(1+g(e^{i\theta}))^{u-v} d\theta.
\end{equation}
Note that
\begin{equation}\label{5.15'}
G_{ab}^\ast (z,w)(1+g(z))^c(1+g(w))^{-d}=G_{a+c,b+d}^\ast (z,w).
\end{equation}
Fix $\epsilon >0$ and let $f^\epsilon (x)=f(x)e^{-\epsilon x}$. Then
$f^\epsilon$ is in $L^1(\mathbb{R})$ and we have
\begin{equation}\label{5.16}
f_N(x)^m=\frac 1{(2\pi )^m} \int_{\mathbb{R}^m} F_m^\epsilon (\lambda)
e^{i\xi_m(\lambda) (x-cN)} d^m\lambda,
\end{equation}
where $F_m^\epsilon (\lambda)=\hat{f}^\epsilon(\lambda_1)\dots 
\hat{f}^\epsilon(\lambda_m)$, $c=2\alpha (1+\alpha)^{-1}$ and
$\xi_m(\lambda) =(\lambda_1+\dots+\lambda_m-im\epsilon)/dN^{1/3}$ with
$d$ given by (\ref{4.6}).
Integration by parts gives
\begin{equation}\label{5.16'}
(1+g(z))^{u-v}-1=(u-v)g(z)+\frac{u-v}{d^4N^{4/3}} R_1(z)+
\frac{(u-v)^2}{d^4N^{4/3}}R_2(z),
\end{equation}
where
\begin{equation}
R_1(z)=d^4N^{4/3}g(z)^2 \int_0^1 \frac{1-t}{(1+tg(z))^2} dt
\notag
\end{equation}
and 
\begin{equation}
R_2(z)=d^4N^{4/3}(\log (1+g(z)))^2 \int_0^1 (1-t)(1+g(z))^{t(u-v)}dt.
\notag
\end{equation}
Let $\Delta f(x)=f(x+1)-f(x)$ be the usual finite difference
operator. We have the following formula
\begin{align}\label{5.17}
&\sum_{x\in\mathbb{Z}} \phi(x,y)\left( \frac{w^y}{z^x}+\frac{w^x}{z^y}
\right) f_N(x)^m=f_N(y)^m \frac{w^y}{z^y} [(1+g(z))^{u-v}
+(1+g(w))^{u-v}]
\\
&-\frac {\alpha (u-v)}{(1-\alpha)^2} \left[\frac{w^y}{z^{y+1}}\Delta
  f_N^m(y) -\frac{w^{y-1}}{z^{y}}\Delta f_N^m(y-1)+
\frac{w^{y+1}}{z^{y}}\Delta f_N^m(y)\right.
\notag\\
&\left.+\frac{w^y}{z^{y-1}}\Delta
  f_N^m(y-1)\right]
+\frac 1{(2\pi)^m} \int_{\mathbb{R}^m} d^m\lambda
F_m^\epsilon(\lambda) e^{i\xi_m(\lambda)(y-cN)}\notag\\
&\left\{ \frac{u-v}{d^4N^{4/3}}\left[ R_1(z
e^{-i\xi_m(\lambda)})-R_1(z)
R_1(w e^{i\xi_m(\lambda)})-R_1(w) \right]\right.\notag\\
&\left.+\frac{(u-v)^2}{d^4N^{4/3}}\left[R_2(z
e^{-i\xi_m(\lambda)})-R_2(z)+R_2(w
e^{i\xi_m(\lambda)})-R_2(w)\right]\right\}
\notag
\end{align}
for $|w|=\exp(-m\epsilon/dN^{1/3})=r_1$, $|z|=r_2=1/r_1$.
To prove this, introduce the formula (\ref{5.151}) for $\phi$ and the
formula (\ref{5.16}) for $f_N^m$ into the left hand side of
(\ref{5.17}) and use 
\begin{equation}
\sum_{x\in\mathbb{Z}} (e^{-i\theta}r_1 e^{i\phi}e^{i\xi_m})^x=
\delta_0(\theta- \frac 1{dN^{1/3}}(\lambda_1+\dots+\lambda_m)-\phi),
\notag
\end{equation}
where $\delta_0$ is the Dirac $\delta$-function, to carry out the
$x$-summation. This gives
\begin{align}
&\sum_{x\in\mathbb{Z}} \phi(x,y)\left( \frac{w^y}{z^x}+\frac{w^x}{z^y}
\right) f_N(x)^m
=\frac 1{(2\pi)^m} \int_{\mathbb{R}^m} d^m\lambda
F_m^\epsilon(\lambda) e^{i\xi_m(\lambda)(y-cN)}\frac{w^y}{z^y}
\\&
\times[(1+g(ze^{-i\xi_m(\lambda)}))^{u-v}+(1+g(we^{i\xi_m(\lambda)}))^{u-v}] 
\notag\\
&=f_N(y)^m \frac{w^y}{z^y} [(1+g(z))^{u-v}
+(1+g(w))^{u-v}]
\notag\\
&+\frac 1{(2\pi)^m} \int_{\mathbb{R}^m} d^m\lambda
F_m^\epsilon(\lambda) e^{i\xi_m(\lambda)(y-cN)}
\frac{w^y}{z^y}
[(1+g(ze^{-i\xi_m(\lambda)}))^{u-v}-(1+g(z))^{u-v}
\notag\\
&+(1+g(we^{i\xi_m(\lambda)}))^{u-v}-(1+g(w))^{u-v}].
\notag
\end{align}
In the last expression we use (\ref{5.16'}) and the explicit form
(\ref{5.141}) of $g$ to obtain the right hand side of (\ref{5.17}). We
will call the first part of the right hand side of (\ref{5.17}),
$f_N(y)^m \frac{w^y}{z^y} [(1+g(z))^{u-v}+(1+g(w))^{u-v}]$, the {\it
  contraction term}, which is the main contribution. The second part
is called the finite difference term.

We can now insert the integral formula (\ref{5.150}) into (\ref{5.14})
and use (\ref{5.17}). The contraction term from the first sum in
(\ref{5.14}) will then exactly cancel the second sum. Here we use
(\ref{5.15'}). What remains is
\begin{equation}\label{5.22}
\frac {u-v}{d^4N^{4/3}} S_0 +\frac {u-v}{d^4N^{4/3}} S_1
+\frac {(u-v)^2}{d^4N^{4/3}} S_2,
\end{equation}
where $S_0$ is the part coming from the finite differences, $S_1$ is
the part coming from terms involving $R_1$ and $S_2$ from the terms
involving $R_2$. After some computation we find
\begin{align}\label{5.23}
S_0&=-\frac{\alpha}{(1-\alpha)^2}\frac 1{(2\pi i)^2}\int_{\gamma_{r_2}}
dz\int_{\gamma_{r_1}}dw G^\ast_{vu}(z,w) (\frac 1z+w)\\&
\times\sum_{y\in\mathbb{Z}} \frac {w^y}{z^y}(dN^{1/3}(f_N(y+1)-f_N(y)))^4.
\notag
\end{align}
Also,
\begin{align}\label{5.231}
S_i&=\frac 1{(2\pi)^3} \int_{\mathbb{R}^3} d^3\lambda F_3^\epsilon
(\lambda) \left(\sum_{y\in\mathbb{Z}} f_N(y)e^{i\xi_3(\lambda)(y-cN)}
  \frac {w^y}{z^y}\right)
\\
&\times\frac 1{(2\pi i)^2}\int_{\gamma_{r_2}}
dz\int_{\gamma_{r_1}}dw G^\ast_{vu}(z,w)\frac {w^y}{z^y}\left\{
  6[h(z,w;\xi_2(\lambda)) \right.
\notag\\
&\left.-h(z,w;0)]-4[h(z,w;\xi_1(\lambda))-h(z,w;0)]-
4[h(z,w;\xi_3(\lambda))-h(z,w;0)]\right\},\notag
\end{align}
$i=1,2$. In order to restrict the $y$-summation so that
$(y-cN)/dN^{1/3}$ ranges over a compact interval we make a summation
by parts in (\ref{5.231}). Recall that we assume that $f(y)$ is a
constant for large $y$. If we let $a=\exp(i\xi_3(\lambda))w/z$ and use
$a^y=(1-a)^{-1}(a^y-a^{y+1})$ in the $y$-sum in (\ref{5.231}) a
summation by parts gives
\begin{equation}
\frac 1{1-\exp(i\xi_3(\lambda))w/z}\sum_{y\in\mathbb{Z}}
(f_N(y)-f_N(y-1))e^{i\xi_3(\lambda)y}\frac {w^y}{z^y}.
\notag
\end{equation} 
Hence, for $i=1,2$,
\begin{align}\label{5.25}
S_i&=
\sum_{y\in\mathbb{Z}} (f_N(y)-f_N(y-1))
\frac 1{(2\pi)^3} \int_{\mathbb{R}^3} d^3\lambda F_3^\epsilon
(\lambda)e^{i\xi_3(\lambda)(y-cN)}
\\
&\frac 1{(2\pi i)^2}\int_{\gamma_{r_2}}
dz\int_{\gamma_{r_1}}dw G^\ast_{vu}(z,w)\frac {w^y}{z^y}\left\{
  6[h(z,w;\xi_2(\lambda)) \right.
\notag\\
&-h(z,w;0)]-4[h(z,w;\xi_1(\lambda))-h(z,w;0)]-\notag\\&\left.
4[h(z,w;\xi_3(\lambda))-h(z,w;0)]\right\}
\frac 1{1-\exp(i\xi_3(\lambda))w/z}
\notag
\end{align}
The expressions $S_i$ will be estimated using the types of estimates
derived in sect. 4.

Write $u=(1+\alpha)(1-\alpha)^{-1}d^{-1}N^{2/3}\tau$,
$v=(1+\alpha)(1-\alpha)^{-1}d^{-1}N^{2/3}\tau'$ and
$y=2\alpha(1-\alpha)^{-1}N (\xi-\tau^2)dN^{2/3}$. To estimate
(\ref{5.23}) we can now use our results from section 4. We use
\begin{align}\label{5.25'}
z(t)&=p_c(\beta)+\frac {\eta}{dN^{1/3}}-\frac {it}{dN^{1/3}}=
1+\frac {\tau+\eta-it}{dN^{1/3}}+\dots
\\
w(t)^{-1}&=p_c(\beta')+\frac {\eta'}{dN^{1/3}}-\frac {is}{dN^{1/3}}=
1+\frac {-\tau'+\eta'-is}{dN^{1/3}}+\dots
\notag
\end{align}
as parametrizations of the integrals as before. Using the same
estimates as in section 4 we can restrict the integration to
$|t|,|s|\le N^\gamma$ with an error $\le C\exp(-cN^{2\gamma})$ with
some $c>0$. Since $v-u\ge 1$, and hence $t_v-t_u\ge c/N^{2/3}$, and
furthermore $|t_u|\le\log N$, we can incorporate the error term into
the right hand side of (\ref{5.2}). The integral in (\ref{5.23}) can
then be estimated using (\ref{4.14}). Note that, by our assumptions on
$f$, the number of $y$-terms $\neq 0$ is $\le CN^{2/3}$ and we get a
compensating factor $1/N^{2/3}$ from the parametrizations; see also
(\ref{4.18}). The numbers $\eta,\eta'$ are chosen so that
$\eta+\eta'\ge 2$. Note that, since we assume $\tau'-\tau\le 1$, the
condition (\ref{4.15}) is satisfied. We find
\begin{equation}\label{5.26}
|S_0|\le C(f,\alpha)\frac 1{N^{2/3}}\sum_y e^{(\tau^3-{\tau'}^3)/3 
+\xi(\tau'-\tau)+(\eta^3+{\eta'}^3)/3-\xi(\eta+\eta')},
\end{equation}
where the $y$-summation is over all $y\in\mathbb{Z}$ such that
$(y-cN)/dN^{1/3}\in [K_1-1,K_2+1]$. 

Consider now $S_i$. Write $\tilde{z}=z\exp(-m\epsilon/dN^{1/3})$,
$\tilde{\lambda}=(\lambda_1+\dots+\lambda_m)/dN^{1/3}$. Then,
\begin{equation}
g(z
e^{-i\xi_m(\lambda)})=-\frac{\alpha}{(1-\alpha)^2}[(\tilde{z}-1)^2\frac
1{\tilde{z}}e^{i\tilde{\lambda}} -2i(\tilde{z}-1)\sin\tilde{\lambda}
+2(\cos\tilde{\lambda}-1)].
\notag
\end{equation}
Thus,

\begin{equation}
1+g(z
e^{-i\xi_m(\lambda)})=1+\frac{2\alpha
  t}{(1-\alpha)^2}(1-\cos\tilde{\lambda}) -\frac{\alpha
  t}{(1-\alpha)^2}[\frac{(\tilde{z}-1)^2}
 {\tilde{z}}e^{i\tilde{\lambda}}-2i(\tilde{z}-1)\sin\tilde{\lambda}].
\notag
\end{equation}
The last term is small for large $N$ and the second is $\ge 1$. Hence,
\begin{equation}\label{5.27}
\frac 1{|1+tg(z(t)e^{-i\xi_m(\lambda)})|}\le 2
\end{equation}
for $|t|\le N^\gamma$ and $N$ sufficiently large. Consequently, there
is a constant $c_1(\alpha)$ depending only on $\alpha$ such that 
\begin{equation}\label{5.28}
|g(ze^{-i\xi_m(\lambda)})|\le
 c_1(\alpha)(\tilde{\lambda}^2+|\tilde{z}-1|^2) 
\end{equation}
and
\begin{equation}\label{5.29}
|\log(1+tg(ze^{-i\xi_m(\lambda)}))|\le
 c_1(\alpha)(\tilde{\lambda}^2+|\tilde{z}-1|^2).
\end{equation}
We can also write
\begin{align}
&1+tg(ze^{-i\xi_m(\lambda)})=(1+\frac
{2\alpha}{(1-\alpha)^2}(1-\cos\tilde{\lambda}))\\ 
&\times\left[1-\frac {\alpha}{(1-\alpha)^2}\frac {(\tilde{z}-1)^2
e^{i\tilde{\lambda}}/\tilde{z}-2i(\tilde{z}-1)\sin\tilde{\lambda}}{
2\alpha (1-\alpha)^{-2}(1-\cos\tilde{\lambda})} \right].
\notag
\end{align}
By periodicity it is enough to consider
$|\tilde{\lambda}|\le\pi$. Estimating the cosine and sine functions we
see that there are constants $c_2(\alpha)$ and $c_3(\alpha)$ such that
\begin{equation}
|1+tg(ze^{-i\xi_m(\lambda)})|
\ge \exp(c_2(\alpha)\tilde{\lambda}^2-c_3(\alpha)(|\tilde{z}-1|^2+
|\tilde{z}-1||\tilde{\lambda}|))
\notag
\end{equation}
and hence,
\begin{equation}
|1+tg(ze^{-i\xi_m(\lambda)})|^{t(u-v)}
\le \exp(-c_2(\alpha)\tilde{\lambda}^2+c_3(\alpha)(|\tilde{z}-1|^2+
|\tilde{z}-1||\tilde{\lambda}|)).
\notag
\end{equation}
Estimating the quadratic polynomial in $\tilde{\lambda}$ we obtain
\begin{equation}
|1+tg(ze^{-i\xi_m(\lambda)})|^{t(u-v)}
\le \exp(t(\tau'-\tau)c_4(\alpha)N^{2/3}|\tilde{z}-1|^2).
\notag
\end{equation}
Now,
$z(t)\exp(-m\epsilon/dN^{1/3})=1+(\tau+\eta-m\epsilon-it)/dN^{1/3}+\dots$,
and we obtain an estimate
\begin{equation}\label{5.30}
|1+tg(z(t)e^{-i\xi_m(\lambda)})|^{t(u-v)}
\le \exp(c_5(\alpha)[(\tau+\eta-m\epsilon)^2+t^2]).
\end{equation}
A computation shows that
\begin{equation}\label{5.31}
\frac 1{|1-e^{i\xi_3(\lambda)w/z}|}\le c_6
\end{equation}
if $\tau-\tau'+\eta+\eta'>\epsilon>0$. Since $\tau'\tau\le 1$ and we
take $\eta+\eta'\ge 2$ we see that we can take $\epsilon=1/2$ for
example. Furthermore, since $(y-cN)/dN^{1/3}$ is bounded for the $y$:s
that contribute to the sum,
\begin{equation}\label{5.32}
|e^{i\xi_3(\lambda)(y-cN)}|\le c_7.
\end{equation}
We can now again estimate as in sect. 4 and use (\ref{4.14}). This
results in an estimate
\begin{align}\label{5.33}
&|S_i|\le
c_8(t,\alpha)\left(\int_{\mathbb{R}}(1+\lambda^2)|\hat{f}^\epsilon
  (\lambda)| d\lambda\right)^3e^{c_5(\alpha)(\tau+\eta-m\epsilon)^2}\\
&\times\frac 1{N^{1/3}}\sum_y
e^{(\tau^3-{\tau'}^3)/3+\xi(\tau'-\tau)+(\eta^3+{\eta'}^3)/3-\xi(\eta+\eta')} 
\notag\\
&\times\int_{\mathbb{R}}e^{(c_5(\alpha)-\eta/2)t^2} dt
\int_{\mathbb{R}}e^{(c_5(\alpha)-\eta'/2)s^2} ds.\notag
\end{align}
We pick $\eta,\eta'\ge 3c_5(\alpha)$. Recall that
$\xi=(y-cN)/dN^{1/3}+\tau^2$. Let $\eta=\max (|\tau|,3c_5(\alpha),1)$
and $\eta'=\max (|\tau'|,3c_5(\alpha),1)$. It follows from
(\ref{5.26}) and (\ref{5.33}) that 
\begin{equation}\label{5.34}
|S_i|\le c_9(f,\alpha),
\end{equation}
$i=1,2$ if $|\tau|,|\tau'|$ are small. If $|\tau|$ and $|\tau'|$ are
large, say $\tau,\tau'\gg 1$, then $\eta=\tau$ and $\eta'=\tau'$,
$0\le\tau'-\tau\le 1$, and we get from (\ref{5.26}) and (\ref{5.33})
that
\begin{equation}\label{5.35}
|S_i|\le c_{10}(f,\alpha)e^{-\tau^3}.
\end{equation}
Inserting these estimates into (\ref{5.22}) and using $v-u\ge 1$, we
obtain an estimate of (\ref{5.14}) of the type we have in the right
hand side of (\ref{5.2}).

Consider the expression
\begin{equation}\label{5.36}
\sum_{x_2,x_4\in\mathbb{Z}} \phi(x_2,x_4)f_N(x_2)f_N(x_4)
\tilde{K}_{cd}^{ab}
\left(\begin{matrix} x_4 & x_3 \\ x_2 & x_1
  \end{matrix}\right).
\end{equation}
In our computations with the kernel $\tilde{K}$ given by (\ref{5.150})
we will leave out the complex integrations. Thus
\begin{align}
&\tilde{K}_{cd}^{ab}
\left(\begin{matrix} x_4 & x_3 \\ x_2 & x_1
  \end{matrix}\right)=
\left|\begin{matrix} G^\ast_{ac}(z_1,w_1)\frac{w_1^{x_2}}{z_1^{x_4}} &
G^\ast_{ad}(z_1,w_1)\frac{w_1^{x_1}}{z_1^{x_4}} \\
G^\ast_{bc}(z_2,w_2)\frac{w_2^{x_2}}{z_2^{x_3}} &
G^\ast_{bd}(z_2,w_2)\frac{w_2^{x_1}}{z_2^{x_3}}
\end{matrix} \right|
\notag\\
&=G^\ast_{ac}(z_1,w_1)G^\ast_{bd}(z_2,w_2)\frac{w_1^{x_2}}{z_1^{x_4}}
\frac{w_2^{x_1}}{z_2^{x_3}}-
G^\ast_{ad}(z_1,w_1)G^\ast_{bc}(z_2,w_2)\frac{w_1^{x_1}}{z_1^{x_4}}
\frac{w_2^{x_2}}{z_2^{x_3}}.
\notag
\end{align}
We are led to the symmetrized expression
\begin{equation}\label{5.36'}
\frac 12\sum_{x_2,x_4\in\mathbb{Z}}\phi(x_2,x_4)f_N(x_2)f_N(x_4)
\left[\frac{w_1^{x_2}}{z_1^{x_4}}+ 
\frac{w_1^{x_4}}{z_1^{x_2}}\right].
\end{equation}
Perform the $x_4$-summation first and use the formula
(\ref{5.17}). The parts containing $R_1$ and $R_2$ can be estimated in
the same way as above. We will only discuss the contraction and
finite-difference parts. The contraction part of (\ref{5.36'}) is
\begin{equation}\label{5.37}
\frac 12\sum_{x_2\in\mathbb{Z}}f_N(x_2)^2
\frac{w_1^{x_2}}{z_1^{x_2}}
[(1+g(z_1))^{u-v}+(1+g(w_1))^{u-v}]
\end{equation}
and hence the contraction part of (\ref{5.36}) is
\begin{align}
\frac 12\sum_{x_2\in\mathbb{Z}}f_N(x_2)^2\{ &G^\ast_{a+u-v,c}(z_1,w_1)
G^\ast_{bd}(z_2,w_2)\frac{w_1^{x_2}}{z_1^{x_2}}
\frac{w_2^{x_1}}{z_2^{x_3}}\notag\\
&+G^\ast_{a,c+v-u}(z_1,w_1)
G^\ast_{bd}(z_2,w_2)\frac{w_1^{x_2}}{z_1^{x_2}}
\frac{w_2^{x_1}}{z_2^{x_3}}\notag\\
&-G^\ast_{a+u-v,d}(z_1,w_1)
G^\ast_{bc}(z_2,w_2)\frac{w_1^{x_1}}{z_1^{x_2}}
\frac{w_2^{x_2}}{z_2^{x_3}}\notag\\
&-G^\ast_{ad}(z_1,w_1)
G^\ast_{b,c+v-u}(z_2,w_2)\frac{w_1^{x_1}}{z_1^{x_2}}
\frac{w_2^{x_2}}{z_2^{x_3}}\}.
\notag
\end{align}
Here we have also used (\ref{5.15'}). Performing the complex
integrations we obtain
\begin{align}\label{5.38}
&\frac 12\sum_{x_2\in\mathbb{Z}}f_N(x_2)^2\{
\tilde{K}_{a+u-v,c}(x_2,x_2)\tilde{K}_{bd}(x_3,x_1)+
\tilde{K}_{a,c+v-u}(x_2,x_2)\tilde{K}_{bd}(x_3,x_1)\\
&-\tilde{K}_{a+u-v,c}(x_2,x_1)\tilde{K}_{bc}(x_3,x_2)
-\tilde{K}_{a,d}(x_2,x_1)\tilde{K}_{b,c+v-u}(x_3,x_2).\notag
\end{align}
The finite difference part of (\ref{5.36}) is
\begin{align}\label{5.39}
&\frac{\alpha}{2(1-\alpha)^2}\sum_{x_2\in\mathbb{Z}}\Delta f_N(x_2)^2\{
G^\ast_{ac}(z_1,w_1)G^\ast_{bd}(z_2,w_2)\frac{w_1^{x_2}}{z_1^{x_2}}
\frac{w_2^{x_1}}{z_2^{x_3}} (\frac 1{z_1}+w_1)\\
&-G^\ast_{ad}(z_1,w_1)G^\ast_{bc}(z_2,w_2)\frac{w_2^{x_2}}{z_1^{x_2}}
\frac{w_1^{x_1}}{z_2^{x_3}} (\frac 1{z_1}+w_2)\}.
\notag
\end{align}
The double $\phi$-term in (\ref{5.5}), $\Sigma_b$ in (\ref{5.8}) and
(\ref{5.11}) combined give
\begin{align}\label{5.40}
&12\sum_x \phi(x_2,x_4)\phi(x_1,x_3)\tilde{K}_{uu}^{vv}
\left(\begin{matrix} x_4 & x_3 \\ x_2 & x_1
  \end{matrix}\right)f_N(x_1)f_N(x_2)f_N(x_3)f_N(x_4)\\
&-12\sum_x \phi(x_1,x_3)\left[\tilde{K}_{uu}^{uv}
\left(\begin{matrix} x_2 & x_3 \\ x_2 & x_1
  \end{matrix}\right)+\tilde{K}_{vu}^{vv}
\left(\begin{matrix} x_2 & x_3 \\ x_2 & x_1
  \end{matrix}\right)\right]
f_N(x_1)f_N(x_2)^2f_N(x_3)\notag\\
&+3\sum_x \left[\tilde{K}_{uu}^{uu}
\left(\begin{matrix} x_1 & x_2 \\ x_1 & x_2  \end{matrix}\right)
+2\tilde{K}_{uv}^{uv}
\left(\begin{matrix} x_1 & x_2 \\ x_1 & x_2  \end{matrix}\right)
+\tilde{K}_{vv}^{vv}
\left(\begin{matrix} x_1 & x_2 \\ x_1 & x_2  \end{matrix}\right)\right]
f_N(x_1)^2f_N(x_2)^2\notag\\
&\doteq A_1+A_2+A_3.
\notag
\end{align}
Consider the $x_4$-summation in $A_1$. The contraction part is, by
(\ref{5.38}),
\begin{equation}
6\sum_x f_N(x_1)f_N(x_2)^2f_N(x_3)\phi(x_1,x_2)
\left[\tilde{K}_{uu}^{uv}
\left(\begin{matrix} x_2 & x_3 \\ x_2 & x_1
  \end{matrix}\right)+\tilde{K}_{vu}^{vv}
\left(\begin{matrix} x_2 & x_3 \\ x_2 & x_1
  \end{matrix}\right)\right],
\notag
\end{equation}
which is exactly $-\frac12 A_2$. Hence what remains of $A_2$ is
\begin{equation}\label{5.41}
-6\sum_x \phi(x_1,x_3)\left[\tilde{K}_{uu}^{uv}
\left(\begin{matrix} x_2 & x_3 \\ x_2 & x_1
  \end{matrix}\right)+\tilde{K}_{vu}^{vv}+
\left(\begin{matrix} x_2 & x_3 \\ x_2 & x_1
  \end{matrix}\right)\right]
f_N(x_1)f_N(x_2)^2f_N(x_3).
\end{equation}
We have
\begin{equation}
\sum_{x_1,x_3\in\mathbb{Z}} f_N(x_1)f_N(x_3)\tilde{K}_{uu}^{uv}
\left(\begin{matrix} x_2 & x_3 \\ x_2 & x_1
  \end{matrix}\right)=
\sum_{x_2,x_4\in\mathbb{Z}} f_N(x_2)f_N(x_4)\tilde{K}_{uu}^{vu}
\left(\begin{matrix} x_4 & x_1 \\ x_2 & x_1
  \end{matrix}\right).
\notag
\end{equation}
We can now apply (\ref{5.38}) to compute the contraction part of the
first half of (\ref{5.41}) and get
\begin{equation}\label{5.42}
-6\sum_{x_1,x_2\in\mathbb{Z}}\left[\tilde{K}_{uu}^{uu}
\left(\begin{matrix} x_1 & x_2 \\ x_1 & x_2
  \end{matrix}\right)+\tilde{K}_{uv}^{uv}+
\left(\begin{matrix} x_1 & x_2 \\ x_1 & x_2
  \end{matrix}\right)\right]
f_N(x_1)^2f_N(x_2)^2.
\end{equation}
Similarly the second half of (\ref{5.41}) has the contraction part
\begin{equation}\label{5.43}
-6\sum_{x_1,x_2\in\mathbb{Z}}\left[\tilde{K}_{uv}^{uv}
\left(\begin{matrix} x_1 & x_2 \\ x_1 & x_2
  \end{matrix}\right)+\tilde{K}_{vv}^{vv}+
\left(\begin{matrix} x_1 & x_2 \\ x_1 & x_2
  \end{matrix}\right)\right]
f_N(x_1)^2f_N(x_2)^2.
\end{equation}
Since the contraction part of (\ref{5.41}) equals (\ref{5.42}) plus
(\ref{5.43}) we see that this exactly cancels $A_3$. It remains to
consider the finite difference parts.

From $A_1$ we get a finite difference part
\begin{align}\label{5.44}
&\frac{6\alpha}{(1-\alpha)^2}\sum_{x}\left[ G^\ast_{vu}(z_1,w_1)
G^\ast_{vu}(z_2,w_2) \frac{w_1^{x_2}}{z_1^{x_2}}
\frac{w_2^{x_1}}{z_2^{x_3}}(\frac 1{z_1}+w_1)\right.\\
&\left.-
G^\ast_{vu}(z_1,w_1)
G^\ast_{vu}(z_2,w_2) \frac{w_2^{x_2}}{z_1^{x_2}}
\frac{w_1^{x_1}}{z_2^{x_3}}(\frac 1{z_1}+w_2)\right]
\Delta f_N(x_2)^2\phi(x_1,x_3)f_N(x_1)f_N(x_3).
\notag
\end{align}
We also need the finite difference part of (\ref{5.41}). These finite
difference parts should be cancelled by the contraction part of
(\ref{5.44}). The contraction part of (\ref{5.41}) is
\begin{align}\label{5.45}
&-\frac{3\alpha}{(1-\alpha)^2}\sum_{x}\left[ G^\ast_{vu}(z_1,w_1)
G^\ast_{uu}(z_2,w_2) \frac{w_1^{x_2}}{z_1^{x_2}}
\frac{w_2^{x_1}}{z_2^{x_1}}(\frac 1{z_1}+w_1)\right.\\
&\left.-
G^\ast_{vu}(z_1,w_1)
G^\ast_{uu}(z_2,w_2) \frac{w_2^{x_2}}{z_1^{x_2}}
\frac{w_1^{x_1}}{z_2^{x_1}}(\frac 1{z_1}+w_2)\right]
\Delta f_N(x_2)^2\phi(x_1,x_3)f_N(x_1)^2
\notag\\
&-\frac{3\alpha}{(1-\alpha)^2}\sum_{x}\left[ G^\ast_{vu}(z_1,w_1)
G^\ast_{vv}(z_2,w_2) \frac{w_1^{x_2}}{z_1^{x_2}}
\frac{w_2^{x_1}}{z_2^{x_1}}(\frac 1{z_1}+w_1)\right.\notag\\
&\left.-
G^\ast_{vv}(z_1,w_1)
G^\ast_{vu}(z_2,w_2) \frac{w_2^{x_2}}{z_1^{x_2}}
\frac{w_1^{x_1}}{z_2^{x_1}}(\frac 1{z_1}+w_2)\right]
\Delta f_N(x_2)^2\phi(x_1,x_3)f_N(x_1)^2.
\notag
\end{align}
In (\ref{5.44}) we have first
\begin{align}
&\sum_x \frac{w_1^{x_2}}{z_1^{x_2}}
\frac{w_2^{x_1}}{z_2^{x_3}}\phi(x_1,x_3)f_N(x_1)f_N(x_3)\Delta
f_N(x_2)^2
\notag\\
&=\frac 12\sum_{x_1}\Delta f_N(x_1)^2\frac{w_1^{x_1}}{z_1^{x_1}}
\left(\sum_{x_2,x_4}\left(\frac{w_2^{x_2}}{z_2^{x_4}}+
\frac{w_2^{x_4}}{z_2^{x_2}}\right)\phi(x_2,x_4)f_N(x_2)f_N(x_4)
\right).
\notag
\end{align}
This gives the contraction part
\begin{equation}\label{5.46}
\frac 12\sum_{x_1,x_2}\frac{w_1^{x_1}}{z_1^{x_1}}
\frac{w_2^{x_2}}{z_2^{x_2}}f_N(x_2)^2\Delta f_N(x_1)^2
[(1+g(z_2))^{u-v}+(1+g(w_2))^{u-v}].
\end{equation}
From the other half of (\ref{5.44}) we get similarly the contraction
part
\begin{equation}\label{5.47}
\frac 12\sum_{x_1,x_2}\frac{w_2^{x_1}}{z_1^{x_1}}
\frac{w_1^{x_2}}{z_2^{x_2}}f_N(x_2)^2\Delta f_N(x_1)^2
[(1+g(z_2))^{u-v}+(1+g(w_1))^{u-v}].
\end{equation}
By (\ref{5.46}) and (\ref{5.47}) the contraction part of (\ref{5.44})
is
\begin{align}
&\frac{3\alpha}{(1-\alpha)^2}\sum_{x}\frac{w_1^{x_1}}{z_1^{x_1}}
\frac{w_2^{x_2}}{z_2^{x_2}}\Delta f_N(x_1)^2f_N(x_2)^2
(\frac 1{z_1}+w_1) 
\notag\\
&\times[G^\ast_{vu}(z_1,w_1)G^\ast_{uu}(z_2,w_2)
+G^\ast_{vu}(z_1,w_1)G^\ast_{vv}(z_2,w_2)]\notag\\
-&\frac{3\alpha}{(1-\alpha)^2}\sum_{x}\frac{w_2^{x_1}}{z_1^{x_1}}
\frac{w_1^{x_2}}{z_2^{x_2}}\Delta f_N(x_1)^2f_N(x_2)^2
(\frac 1{z_1}+w_2) 
\notag\\
&\times[G^\ast_{vu}(z_1,w_1)G^\ast_{uu}(z_2,w_2)
+G^\ast_{vv}(z_1,w_1)G^\ast_{vu}(z_2,w_2)],
\notag
\end{align}
which exactly cancels (\ref{5.45}). The finite difference part of
(\ref{5.44}) is handled in the same way as (\ref{5.23}).

We will end with some brief comments about remaining estimates.
By (\ref{5.150}) we have for example
\begin{align}
\Delta K_{uv}(x,y)&=\frac 1{(2\pi i)^2} \int_{\gamma_{r_2}} dz
\int_{\gamma_{r_1}} dw \left(\frac {1-\alpha/z}{1-\alpha z} \right)^N 
\left(\frac {1-\alpha w}{1-\alpha/w} \right)^N \frac {w^y}{z^x} 
\\&\times\frac 1{w(z-w)} [(1+g(w))^{v-u}-1](1+g(z))^u(1+g(w))^u.
\notag
\end{align}
Here we can expand $(1+g(w))^{v-u}-1$ as in (\ref{5.16'}) and then
estimate in the same way as we did for the $R_1-$ and $R_2-$ terms
above. In this way we will see that the $\Delta K^2$-terms give
contributions of the right type. We get a similar integral expression
for $\Delta K_{uv}+\Delta K_{vu}-\Delta K_{vv}$ involving
\begin{equation}
(1+g(z))^v(1+g(w))^{-u}[(1+g(z))^{u-v}-1][(1+g(w))^{v-u}-1]
\notag
\end{equation}
and we proceed similarly.

\end{proof}

\subsection{Weak convergence}
Consider $H_N(f,t_j)$ as defined by (\ref{5.1}). The next lemma is a
standard consequence of lemma \ref{L5.1}.

\begin{lemma}\label{L5.2}
Under the same assumptions as in lemma \ref{L5.1},
\begin{equation}\label{5.201}
\mathbb{P}[\max_{j=u,\dots,v} |H_N(f,t_j)-H_N(f,t_u)|\ge\lambda] \le
C(f,\alpha)\lambda^{-4} e^{-|t_u|^3}|t_u-t_v|^2.
\end{equation}
\end{lemma}

\begin{proof}
Let
\begin{equation}
\eta_i=H_N(f,t_{u+i})-H_N(f,t_{u+i-1}),
\notag
\end{equation}
$T_m=\sum_{i=1}^m \eta_i$ and $T_0=0$, so that
$T_j-T_i=H_N(f,t_j)-H_N(f,t_i)$. It follows from (\ref{5.2}) and
Chebyshev's inequality that
\begin{equation}
\mathbb{P}[|T_i-T_j|\ge\lambda] \le C(f,\alpha)\lambda^{-4}
e^{-|t_u|^3} \left(\sum_{i<\ell\le j} u_\ell\right)^2
\notag
\end{equation}
for $u\le i<j\le v$, where $u_\ell=t_{\ell+u}-t_{\ell-1+u}$. This
implies (\ref{5.201}) according to theorem 12.2 in \cite{Bill}.
\end{proof}

Fix $l>0$ and consider rescaled top height curve $H_{N,0}(t)$ for
$|t|\le L$ and its modulus of
continuity, $0<\delta\le 1$,
\begin{equation}
w_N(\delta)=\sup_{|t|,|s|\le T,|s-t|\le\delta}
|H_{N,0}(t)-H_{N,0}(s)|,
\notag
\end{equation}
\begin{lemma}\label{L5.3} Let $w_N$ be defined as above.
Given $\epsilon,\lambda>0$ there is a $\delta>0$ and an integer $N_0$
such that
\begin{equation}
\mathbb{P}[w_N(\delta)\ge\lambda]\le\epsilon
\notag
\end{equation}
if $N\ge N_0$.
\end{lemma}
Together with the convergence of the finite dimensional distributions,
theorem \ref{T4.3}, this proves theorem \ref{T0.4}, \cite{Bill}.
We turn now to the proof of lemma \ref{L5.3}.

\begin{proof} Assume that $\delta^{-1},T\in\mathbb{Z}$ and divide the
  interval $[-T,T]$ into $2m$ parts of length $T/m=\delta$. Write
  $r_j=[j\delta[cN^{2/3}]]$, $c=2\alpha/(1-\alpha)$, so that
  $t_{r_j}\approx j\delta$.

\begin{claim}\label{Cl5.4}
Let $L=T[cN^{2/3}]$ and
$B_M$ is the subset of our probability space where 
$\max_{|j|\le L} |H_{N,0}(t_j)|\le M$.
Then, given $\epsilon>0$, we can choose $M$ so that 
\begin{equation}\label{5.203}
\mathbb{P}[B_M^c]\le\epsilon.
\end{equation}
\end{claim}
We will prove this claim below. We will also need
\begin{claim}\label{Cl5.5}
For any $\lambda>0$ there is a constant $C(M)$ that depends on $M$ but
not on $\lambda$ such that
\begin{equation}\label{5.204}
\mathbb{P}[\max_{r_j\le i\le
  r_{j+1}}|H_{N,0}(t_i)-H_{N,0}(t_{r_j})|\ge \lambda, B_M]\le \frac
  {C(M)}{\lambda}.
\end{equation}
\end{claim}
We will return to the proofs. The proof of both claims are based on
choosing appropriate functions $f$ in lemma \ref{L5.2} and results
about convergence in distribution. Assuming the validity of the two
claims we can prove lemma 5.3. Set
\begin{equation}\label{5.205}
A_j=\{\max_{t_{r_j}\le i\le
  t_{r_{j+1}}}|H_{N,0}(s)-H_{N,0}(t_{r_j})|\ge \lambda/3\},
\end{equation}
so that $\{w_N(\delta)\ge\lambda\}\subseteq\cup_{|j|\le} A_j$.
Choose $M$ so large that $\mathbb{P}[B_M^c]\le\epsilon$, which is
possible by claim \ref{Cl5.4}. Hence
\begin{equation}\label{5.206}
\mathbb{P}[w_N(\delta)\ge \lambda]\le \epsilon +\sum_{|j|\le m}
\mathbb{P}[A_j\cup B_M].
\end{equation}
Now, if the inequality in (\ref{5.205}) holds then
\begin{equation}
\max_{r_j\le i\le
  r_{j+1}}|H_{N,0}(t_i)-H_{N,0}(t_{r_j})|\ge \lambda/9.
\notag
\end{equation}
Consequently, using (\ref{5.204}) and (\ref{5.206}),
\begin{equation}
\mathbb{P}[w_N(\delta)\ge \lambda]\le \epsilon
+(2m+1)C(M)\lambda^{-1}\delta^2 \le \epsilon+2TC(M)\lambda^{-1}\delta.
\notag
\end{equation}
Choose $\delta$ so that $\delta\le\epsilon\lambda/2TC(M)$. Lemma
\ref{L5.3} is proved

Consider now claim \ref{Cl5.4}. Pick a $C^\infty$ function $q$ such
that $0\le g\le 1$ and
\begin{equation}\label{5.206'}
g(x)=
\begin{cases}
1
&\text{if $x\ge 0$}\\
0
&\text{if $x\le -1$}
\end{cases}
\end{equation}
and let $g_M(x)=g(x-M)$. It is not hard to see that if we take $f=g_M$
in (\ref{5.2}) the $C(f,\alpha)$ can be taken to be independent of $M$
(only sup-norms of $f$ and its derivatives enter).
If $H_N(g_M,t_{r_j})>1/4$, then $H_{N,0}(t_{r_j})\ge M-1$ and using
the convergence in distribution to $F_2$ we see that we can choose $M$
so large that
\begin{equation}
\mathbb{P}[H_N(g_M,t_{r_j})>1/4]\le\epsilon^2
\notag
\end{equation}
for $|j|\le m$ and all sufficiently large $N$. 
Let $\omega$ denote a point in our probability space.
Now,
\begin{align}
&\mathbb{P}[\max_{|j|\le L}
H_N(g_M,t_j)>1/2]=\mathbb{P}[H_N(g_M,t_{j(\omega)})>1/2]
\notag\\
&=\sum_{j=-m+1}^m \mathbb{P}[H_N(g_M,t_{j(\omega)})>1/2, t_{r_{j-1}}\le
t_{j(\omega)}\le t_{r_j}] 
\notag\\
&\le \sum_{j=-m+1}^m\mathbb{P}[H_N(g_M,t_{r_jj})>1/4] \notag\\&+
\sum_{j=-m+1}^m\mathbb{P}[ \max_{r_{j-1}\le i\le t_{r_j}}
|H_N(g_M,t_i) -H_N(g_M,t_{r_{j-1}})|\ge 1/4]
\notag\\
&\le 2m\epsilon^2+\sum_{j=-m+1}^m\frac C{(1/4)^4}\delta^2\le
2T\frac{\epsilon^2}{\delta} +C\delta
\notag
\end{align}
by lemma \ref{L5.2}. We can now choose $\delta=\epsilon$.

If $H_{N,0}(t_j)\ge M$, then $H_N(g_m,t_j)\ge 1$ and hence
$
\max_{|j|\le L}H_{N,0}(t_j)\ge M,
$
which implies 
$
\max_{|j|\le L}H_{N}(g_Mt_j)\ge 1/2
$.
It follows that 
\begin{equation}
\mathbb{P}[\max_{|j|\le L}H_{N,0}(t_j)\ge M]\le\epsilon.
\notag
\end{equation}
The case 
$
\max_{|j|\le L}H_{N,0}(t_j)\le -M,
$
is analogous. This proves claim \ref{Cl5.4}.

To prove claim \ref{Cl5.5}  we let $i(\omega)$ be defined by
\begin{equation}
\max_{r_j\le i\le r_{j+1}} |H_{N,0}(t_i)-H_{N,0}(t_{r_j})|=
|H_{N,0}(t_{i(\omega)})-H_{N,0}(t_{r_j})|.
\notag
\end{equation}
Let $I_j=[j\lambda,(j+1)\lambda )$, $j=-K,\dots, K-1$, where
$M=K\lambda$, $K\in\mathbb{Z}^+$. Take a $C^\infty$ function $f$, $0\le
f\le 1$, soch that $f(x)=0$ if $x\le-\lambda$, $f(x)=1$ if $0\le
x\le\lambda$ and $F(x)=0$ if $x\ge\lambda$. Set
\begin{equation}\label{5.206''}
f_j(x)=f(x-\lambda j).
\end{equation}
Suppose first that $H_{N,0}(t_{i(\omega)})\le
H_{N,0}(t_{r_j})-2\lambda$ and that $\omega\in B_M$. Then there is a
$k(\omega)$ such that $H_{N,0}(t_{r_j})\in I_{k(\omega)}$, and
\begin{equation}
H_{N,0}(t_{i(\omega)})\le (k+1)\lambda-2\lambda=(k-1)\lambda,
\notag
\end{equation}
and consequently $H_N(f_{k(\omega)},t_{i(\omega)})=0$. Since
$H_N(f_{k(\omega)},t_{r_j})\ge 1$, we see that
\begin{equation}
|H_N(f_{k(\omega)},t_{i(\omega)})-H_N(f_{k(\omega)},t_{r_j})|\ge 1.
\notag
\end{equation}
Hence,
\begin{equation}\label{5.207}
\max_{|k|\le m}\max_{r_j\le i\le r_{j+1}}
|H_N(f_k,t_i)-H_N(f_k,t_{r_j})|\ge 1.
\end{equation}
Call this event $F$. If we instead suppose that 
$H_{N,0}(t_{i(\omega)})\ge
H_{N,0}(t_{r_j})+2\lambda$ we can proceed similarly and see that
(\ref{5.207}) still holds. Now,
\begin{align}
\mathbb{P}[F]&\le \mathbb{P}[\cup_{r_j\le i\le r_{j+1}} \{
|H_N(f_k,t_i)-H_N(f_k,t_{r_j})|\ge 1\}]
\notag\\
&\le \sum_{k=-K+1}^K C|t_{r_{j+1}}-t_{r_j}|^2\le \frac
{2MC}{\lambda}\delta^2,
\notag
\end{align}
by lemma \ref{L5.2}.

\end{proof}

\subsection{Transversal fluctuations}
In this section we will prove corollary \ref{C0.5}, proposition
\ref{P0.6} and theorem \ref{T0.8}. Let $T>0$ be fixed and set
\begin{align}
S_N^T(u)&=\sup_{-T\le t\le u} H_{N,0}(t)\notag
\\
S^T(u)&=\sup_{-T\le t\le u} (A(t)-t^2).
\notag
\end{align}
We will write $S_N^T$ for $S_N^T(T)$ and $S^T$ for $S^T(T)$.

\begin{lemma}\label{L5.6}
Given $\epsilon>0$ we can choose $T=T(\epsilon)$  so that
\begin{equation}\label{5.208}
\mathbb{P}[S_N^\infty\neq S_N^T]\le \epsilon
\end{equation}
for all sufficiently large $N$.
\end{lemma}

Note that together with theorem \ref{T0.4} this proves corollary
\ref{C0.5}.

\begin{proof}
Let $g_M$ be defined as above and set $R_j=T+(j-1)\delta$,
$j\ge 1$, where $\delta$ will be specified below. It follows from
lemma \ref{L5.2} that for $R_j\le\log N$,
\begin{equation}\label{5.209}
\mathbb{P}[\sup_{R_j\le t\le R_{j+1}} |H_N(g_M,t)-H_N(g_M,R_j)|\ge
1/2]\le Ce^{-R_j^3}\delta^2.
\end{equation}
where $C$ is independent of $M$. Now,
\begin{align}\label{5.210}
&\mathbb{P}[\sup_{T\le t\le R_L} H_{N,0}(t)\ge M]\le
\mathbb{P}[\sup_{T\le t\le R_L} H_{N}(g_M,t)\ge 1]
\\
&=\mathbb{P}[\max_{1\le j<L}\sup_{R_j\le t\le R_{j+1}}H_{N}(g_M,t)\ge
1]
\notag\\
&\le \sum_{j=1}^{L-1} \mathbb{P}[\sup_{R_j\le t\le
  R_{j+1}}(H_{N}(g_M,t) -H_{N}(g_M,R_j)) +H_{N}(g_M,R_j)\ge 1]
\notag\\
&\le \sum_{j=1}^{L-1} \mathbb{P}[\sup_{R_j\le t\le
  R_{j+1}}(H_{N}(g_M,t) -H_{N}(g_M,R_j))\ge 1/2]
+\sum_{j=1}^{L-1} \mathbb{P}[H_{N}(g_M,R_j)\ge 1/2].
\notag
\end{align}

\begin{claim}\label{Cl5.7}
There is a positive constant $c_1$ such that
\begin{equation}\label{5.211}
\mathbb{P}[H_{N,0}(R)\ge s]\le e^{-c_1(s+R^2)^{3/2}}.
\end{equation}
\end{claim}
\begin{proof}
Let $c=(1+\alpha)(1-\alpha)^{-1}d^{-1}$,
$d^3=\alpha(1+\alpha)(1-\alpha)^{-3}$ and $t_j=j/cN^{2/3}$ as before. We have
\begin{align}
\mathbb{P}[H_{N,0}(R)\ge s]&=\mathbb{P}[\frac 1{dN^{1/3}}(h_0(2j,2N-1)
-\frac {2\alpha}{1-\alpha}N)\ge s]\notag\\
&=\mathbb{P}[G(N+cN^{2/3}R,N-cN^{2/3}R)\ge\frac {2\alpha}{1-\alpha}N
sdN^{1/3}]
\notag 
\end{align}
by (\ref{0.19}) and the definition of $H_{N,0}$. Recall that the
parameter $q$ in the geometric distribution $=\alpha^2$. By Corollary
2.4 in \cite{Jo1} we have, for all $K\ge 1$ and $\gamma\ge 1$,
\begin{equation}\label{x.1}
\mathbb{P}[G([\gamma K],K)>Kt]\le e^{-2KJ(t+1)},
\end{equation}
where the function $J$ satisfies
\begin{equation}\label{x.2}
J((1+\sqrt{q\gamma})^2(1-q)^{-1}+\delta)\ge c_1'\delta^{3/2}
\end{equation}
for $0\le \delta\le 1$; $c_1'$ is a positive constant. We take
$K=N-cN^{2/3}R$, $\gamma=(N+cN^{2/3}R)/K$ and
$t=(2\alpha(1-\alpha)^{-1} N+sdN^{1/3})/K$. Pick $\delta$ so that
\begin{equation}
(1+\sqrt{q\gamma})^2(1-q)^{-1}+\delta=1+t.
\notag
\end{equation}
This gives $\delta=dN^{-2/3}(s+R^2)+O(N^{-1})$ and if we insert this
into (\ref{x.2}), the estimate (\ref{x.1}) gives us exactly what we
want.
\end{proof}

If $H_N(g_M,R_j)\ge 1/2$, then $H_{N,0}(R_j)\ge M-1$ and hence
\begin{align}\label{5.212}
\sum_{j=1}^{L-1} \mathbb{P}[H_{N}(g_M,R_j)\ge 1/2]&\le \frac 1{\delta}
\sum_{j=1}^{L-1} e^{-c(M-1+(T+(j-1)\delta)^2)^{3/2}} \delta
\\
&\le \frac 1{\delta} \int_{T-1}^\infty e^{-c(M-1+x^2)^{3/2}} dx.
\notag
\end{align}
Using (\ref{5.209}) we find
\begin{align}\label{5.213}
&\sum_{j=1}^{L-1} \mathbb{P}[\sup_{R_j\le t\le
  R_{j+1}}(H_{N}(g_M,t) -H_{N}(g_M,R_j))\ge 1/2]
\\
&\le\sum_{j=1}^{L-1} Ce^{-R_j^3}\delta^2 \le C\delta\int_{T-1}^\infty e^{-x^3}dx.
\notag
\end{align}
Inserting (\ref{5.212}) and (\ref{5.213}) into (\ref{5.210}) gives
\begin{equation}\label{5.214}
\mathbb{P}[\sup_{T\le t\le R_L} H_{N,0}(t)\ge M]\le \frac 1{\delta}
\int_{T-1}^\infty e^{-c(M-1+x^2)^{3/2}} dx+
C\delta\int_{T-1}^\infty e^{-x^3}dx
\end{equation}
if $R_L\le\log N$. We can take $\delta=1$.

It follows from (\ref{5.212}) that
\begin{align}\label{5.215}
&\mathbb{P}[\sup_{R_L\le t} H_{N,0}(t)\ge M]\le \sum_{t_u\ge R_L}
\mathbb{P}[H_{N,0}(t_u) \ge M]
\\
&\le\sum_{u\ge cN^{2/3}R_L} e^{-c(M+t_u^2)^{3/2}} \le CNe^{-c(\log N)^3}
<\epsilon/4 
\notag
\end{align}
if $N$ is sufficiently large. We know that $\mathbb{P}[H_{N,0}(0)\le
M]\to F_2(M)$ as $N\to\infty$ and we can choose $M$ so large that the
right hand side of (\ref{5.214}) is $\le \epsilon/4$. Together with
(\ref{5.215}) this gives (using symmetry),
\begin{equation}\label{5.216}
\mathbb{P}[\sup_{|t|\ge T} H_{N,0}(t)\ge M]\le\epsilon.
\end{equation}
If $H_{N,0}(0)>M$ and $\sup_{|t|\ge T} H_{N,0}(t)\ge M$, then
$S_N^\infty=S_N^T$ and consequently
\begin{equation}
\mathbb{P}[S_N^\infty\neq S_N^T]\le \mathbb{P}[H_{N,0}(0)\le M]+
\mathbb{P}[\sup_{|t|\ge T} H_{N,0}(t)>M]\le 2\epsilon
\notag\end{equation}
for all sufficiently large $N$.
\end{proof}

We turn now to the transversal fluctuations and the proof of 
theorem \ref{T0.8}. 
Define
\begin{align}
K_N^T&=\inf\{ u\ge -T\,;\,S_N^T(u)=S_N^T\}\notag\\
K^T&=\inf\{ u\ge -T\,;\,S^T(u)=S^T\},
\notag
\end{align}
which give the leftmost point of maximum in $[-T,T]$ before and after
the limit.
We first prove proposition \ref{P0.6}. 
\begin{proof} ({\it Proposition \ref{P0.6}}).
Note that
\begin{equation}
\{K_N<-T\}\subseteq \{\sup_{t\le -T} H_{N,0}(t)\ge H_{N,0}(0)\}.
\notag
\end{equation}
It follows that
\begin{equation}
\mathbb{P}[K_N<-T]\le \mathbb{P}[ H_{N,0}(t)\ge
M]+\mathbb{P}[H_{N,0}(0)<M\} \le 2\epsilon,
\notag
\end{equation}
by (\ref{5.216}) and the discussion proceeding it. Also,
$\{K_N>T\}\subseteq \{S_N^\infty\neq S_N^T\}$ and we can use lemma
\ref{L5.6}.
\end{proof}

\begin{proof} ({\it Theorem \ref{T0.8}}). It follows from lemma
    \ref{L5.6} that given $\epsilon>0$ we can choose $T$ and $N_0$ so that
\begin{equation}\label{5.218'}
\mathbb{P}[K_N=K_N^T]\ge1-\epsilon
\end{equation}
if $N\ge N_0$. Let $h_T:C(\mathbb{R})\to\mathbb{R}$ be defined by
\begin{equation}
h_T(x)=\inf\{u\ge -T;\sup_{-T\le t\le u} x(t)=\sup_{-T\le t\le T} x(t)\},
\notag
\end{equation}
and let
\begin{equation}
D_{h_T}=\{x\in C(\mathbb{R})\,;\, \text{$h_T$ is discontinuous at $x$}\}
\notag
\end{equation}
It follows from our assumption that $\mathbb{P}[D_{h_T}]=0$, since
$h_T$ is continuous at $x$ unless $x$ has two distinct maximum
points. Since $H_{N,0}$ converges in distribution to $X$ in $C[-T,T]$
it follows that
\begin{equation}\label{5.219}
K_N^T=h_T(H_{N,0})\rightarrow h_T(X)=K_T
\end{equation}
as $N\to\infty$. 

Let $\mathcal{D}_T$ be all points of disconituity for $x\to
\mathbb{P}[K^T\le x]$, $T\in\mathbb{Z}$, and $\mathcal{D}=\cup_{T\ge
  1} \mathcal{D}_T$. We will prove that 
\begin{equation}\label{5.220}
\mathbb{P}[K_N\le x]\to\mathbb{P}[K\le x]
\end{equation}
as $N\to\infty$ for all $x\in\mathbb{R}\setminus \mathcal{D}$,
which implies what we want since $\mathcal{D}$ is countable.
All the results and assumptions that are behind the estimate
(\ref{5.218'}) can also be proved for the limiting Airy process and we
can assume that $N_0$ and $T\in\mathbb{Z}_+$ are chosen so that also
\begin{equation}\label{5.221}
\mathbb{P}[K=K^T]\ge 1-\epsilon
\end{equation}
if $N\ge N_0$. Let $x\in\mathbb{R}\setminus \mathcal{D}$. Then,
\begin{equation}
\mathbb{P}[K_N\le x]=\mathbb{P}[K_N^T\le x, K_N^T=K_N]+
\mathbb{P}[K_N\le x,K_N^T\neq K_N],
\notag
\end{equation}
and similarly for $K_N^T$. Hence,
\begin{equation}
|\mathbb{P}[K_N\le x]-\mathbb{P}[K_N^T\le x]|\le 2\epsilon
\notag
\end{equation}
if $N\ge N_0$. Since $x\in\mathbb{R}\setminus \mathcal{D}$ it follows
from (\ref{5.219}) that we can choose $N_1$ so that
$|\mathbb{P}[K_N^T\le x] -\mathbb{P}[K^T\le x]|\le\epsilon$ if $N\ge N_1$.
It follows from (\ref{5.221}) that 
$|\mathbb{P}[K\le x] -\mathbb{P}[K^T\le x]|\le\epsilon$. Combining the
estimates we see that
\begin{equation}
|\mathbb{P}[K_N\le x]-\mathbb{P}[K\le x]|\le 4\epsilon
\notag
\end{equation}
if $N\ge\max(N_0,N_1)$, which proves (\ref{5.220}).
\end{proof}

\noindent
{\bf Acknowledgement}: I thank Peter Forrester for drawing my attention 
a few years ago
to the relation between the exponents occuring in \cite{FNH}  and
\cite{Jo4}.

\end{document}